
\magnification 1200
\hsize = 14.5cm
\hoffset -0.5cm

\font\Bbb=msbm10
\def\BBB#1{\hbox{\Bbb#1}}
\font\Frak=eufm10
\def\frak#1{{\hbox{\Frak#1}}}

\hyphenation{Bor-cherds}
\hyphenation{pre-print}

\def\Ldg{1.1}
\def\Ldk{1.2}
\def\Ldd{1.3}
\def\pairDK{1.4}

\def\vir{2.1}
\def\comm{2.2}
\def\degr{2.3}
\def\Borc{2.4}
\def\Bora{2.5}
\def\Borb{2.6}
\def\qas{2.7}
\def\qass{2.8}
\def\omd{2.9}
\def\omdeg{2.10}
\def\Dab{2.11}
\def\Da{2.12}
\def\Dka{2.13}
\def\skeww{2.14}
\def\Ytens{2.15}
\def\omtens{2.16}
\def\vla{2.17}
\def\Y{2.18}
\def\LL{2.19}
\def\LI{2.20}   
\def\II{2.21}
\def\LLz{2.22}
\def\LIz{2.23}
\def\IIz{2.24}
\def\isom{2.25}

\def\twi{3.1}
\def\xoe{3.2}
\def\Ye{3.3}
\def\Yw{3.4}
\def\MHyp{3.5}
\def\ggz{3.6}
\def\T{3.7}

\def\dko{4.1}
\def\dka{4.2}
\def\dgg{4.3}
\def\dda{4.4}
\def\ddo{4.5}
\def\ko{4.6}
\def\ka{4.7}
\def\ggg{4.8}
\def\da{4.9}
\def\tdd{4.10}
\def\commm{4.11}
\def\Br{4.12}
\def\omrumu{4.13}
\def\omrumuone{4.14}
\def\omrumuo{4.15}
\def\donda{4.16}
\def\Ib{4.17}
\def\dondo{4.18}
\def\ddos{4.19}
\def\BIV{4.20}

\def\wdoa{5.1}
\def\dwd{5.2}
\def\prj{5.3}
\def\expr{5.4}
\def\embd{5.5}
\def\relV{5.6}
\def\relM{5.7}
\def\relVY{5.8}
\def\cab{5.9}
\def\doo{5.10}
\def\wdaa{5.11}
\def\strng{5.12}

\def\ltens{2.1}
\def\voa{2.2}
\def\sing{2.3}
\def\hvmod{2.4}
\def\omi{2.5}

\def\Hone{3.1}
\def\Htwo{3.2}
\def\omdd{3.3}
\def\EE{3.4}
\def\omE{3.5}

\def\bbs{4.1}
\def\main{4.2}
\def\tormod{4.3}
\def\Ydgg{4.4}
\def\Ydkk{4.5}
\def\Ydda{4.6}
\def\Yddo{4.7}

\def\irre{5.1}
\def\emb{5.2}
\def\Memb{5.3}
\def\presr{5.4}
\def\maind{5.5}
\def\irrd{5.6}
\def\DK{5.7}

\def\L{{\cal L}}
\def\D{{\cal D}}
\def\CC{{\cal C}}
\def\U{{\cal U}}
\def\F{{\cal F}}
\def\HVir{{{\cal H}{\cal V}{\mit ir}}}

\def\d{\partial}
\def\g{{\frak g}}
\def\gl{{\it gl}}
\def\glN{{\it gl}_N}
\def\wgl{{\widehat {gl}_N}}
\def\slN{{{\it sl}_N}}

\def\wsl{{\widehat {sl}_N}}
\def\dg{{\dot \g}}
\def\dh{{\dot {\frak h}}}
\def\wdg{{\widehat \dg}}
\def\td{\tilde d_0}
\def\wda{{{\widehat d}_a}}
\def\wL{{\overline L}}
\def\wcL{{{\overline c}_L}}
\def\wCL{{{\overline C}_L}}

\def\om{\omega}
\def\ot{\otimes}
\def\eru{e^{\r\u}{}}
\def\emu{e^{\m\u}{}}
\def\ermu{e^{(\r+\m)\u}{}}

\def\div{{\rm div}}
\def\gdiv{{\g_\div}}
\def\Ddiv{\D_\div}
\def\tor{{\mit tor}}

\def\Hyp{{\mit Hyp}}
\def\hyp{{\mit Hyp}}

\def\Hei{{{\cal H}{\mit ei}}}
\def\R{{\cal R}}
\def\K{{\cal K}}

\def\C{\BBB C}
\def\o{{\bf 1}}
\def\Z{\BBB Z}

\def\r{{\bf r}}
\def\m{{\bf m}}
\def\u{{\bf u}}
\def\vv{{\bf v}}
\def\t{{\bf t}}

\def\z{\left[ z_1^{-1} \delta \left( {z_2 \over z_1} \right) \right]}

\def\zd{\left[ z_1^{-1}  {\d \over \d z_2} 
\delta \left( {z_2 \over z_1} \right) \right]}

\def\zdd{\left[ z_1^{-1} \left( {\d \over \d z_2} \right)^2
\delta \left( {z_2 \over z_1} \right) \right]}

\def\zddd{\left[ z_1^{-1} \left( {\d \over \d z_2} \right)^3
\delta \left( {z_2 \over z_1} \right) \right]}

\def\dzb{{\d\over\d z_2}}
\def\dzbb{{\left(\d\over\d z_2\right)^2}}

\def\End{{\rm End}}
\def\Id{{\rm Id}}
\def\Ind{{\rm Ind}}
\def\rank{{{\rm rank}\hbox{\hskip 0.1cm}}}

\def\Der{{\rm Der}}
\def\Vir{{{\cal V}{\mit ir}}}
\def\wVir{{\overline \Vir}}

\def\deg{{{\rm deg}\hbox{\hskip 0.1cm}}}        
\def\char{{{\rm char}\hbox{\hskip 0.1cm}}}
\def\tr{{\rm tr}}
\def\S{{\left< S \right>}}

\

\hfill
{\sl To Robert Moody}

\

\

\centerline
{\bf Energy-momentum tensor for the toroidal Lie algebras.}

\centerline{
{\bf Yuly Billig}
\footnote{*}{Research supported by the  Natural Sciences and
Engineering Research Council of Canada.}
}

\centerline{School of Mathematics \& Statistics}
\centerline{Carleton University}
\centerline{1125 Colonel By Drive}
\centerline{Ottawa, Ontario, K1S 5B6, Canada}
\centerline{e-mail: billig@math.carleton.ca}

\

\

\

{\bf Abstract.} We construct vertex operator representations
for the full $(N+1)$-toroidal Lie algebra $\g$. We associate with $\g$ a toroidal 
vertex operator algebra, which is a tensor product of an affine VOA,
a sub-VOA of a hyperbolic lattice VOA, affine $\slN$ VOA and a twisted
Heisenberg-Virasoro VOA. The modules for the toroidal VOA are also modules
for the toroidal Lie algebra $\g$.
We also construct irreducible modules for an important
subalgebra $\gdiv$ of the toroidal Lie algebra that corresponds to the divergence
free vector fields. This subalgebra carries a non-degenerate invariant bilinear
form. The VOA that controls the representation theory of $\gdiv$ is a tensor product
of an affine VOA $V_\wdg(c)$ at level $c$, a sub-VOA of a hyperbolic lattice VOA,
affine $\slN$ VOA and a Virasoro VOA at level $\wcL$ with the following condition on 
the central charges: $2(N+1) + \rank V_\wdg(c) + \wcL = 26$.

\

\

{\bf 0. Introduction.}

\

Toroidal Lie algebras are very natural multi-variable generalizations of 
affine Kac-Moody algebras. The theory of affine Lie algebras is rich and 
beautiful, and has many important applications in physics. By large, 
applications of toroidal Lie algebras in physics are still to be discovered. 
We should mention however the papers [IKUX], [IKU], where the toroidal symmetry
is discussed in the context of a 4-dimensional conformal field theory.
We hope that the development of the representation theory of toroidal Lie algebras
will help to find the proper place for these algebras in physical theories.

The construction of a toroidal Lie algebra is totally parallel to the well-known
construction of an (untwisted) affine Kac-Moody algebra [K1]. One starts with 
a finite-dimensional simple Lie algebra $\dg$ and considers maps from 
an $N+1$-dimensional torus into $\dg$. We may identify the algebra of functions
on a torus with the Laurent polynomial algebra 
$\R = \C[t_0^\pm, t_1^\pm, \ldots, t_N^\pm]$, by taking the Fourier basis, setting
$t_k = e^{ix_k}$. The Lie algebra of the $\dg$-valued maps from a torus will then
become $\C[t_0^\pm, t_1^\pm, \ldots, t_N^\pm] \ot \dg$. When $N = 0$, this yields
the loop algebra. 

Just as in affine case, one builds the universal central extension 
$(\R \ot \dg) \oplus \K$ of $\R \ot \dg$. 
However  when $N \geq 1$, the center $\K$
is infinite-dimensional. The infinite-dimensional center
makes this Lie algebra highly degenerate. One can show, for example, that in 
an irreducible bounded weight module, most of the center should act trivially.
To eliminate this degeneracy, we add the Lie algebra of the vector fields on a torus,
$\D = \Der (\R)$ to $(\R \ot \dg) \oplus \K$. The resulting algebra,
$$\g = (\R \ot \dg) \oplus \K \oplus \D$$
is called the toroidal Lie algebra (see Section 1 for details). 
The action of $\D$ on $\K$ is non-trivial,
making the center of the toroidal Lie algebra $\g$ finite-dimensional. This enlarged
algebra will have a much better representation theory.

A major obstruction to building the representation theory of toroidal Lie algebras
is that these algebras, being $\Z^{N+1}$ graded, do not possess a triangular
decomposition whenever $N \geq 1$. For this reason, the standard construction of 
the highest weight modules fails to work. The modules for the toroidal Lie algebras
built by forcing the highest weight condition are not attractive [BC].

The true representation theory for toroidal Lie algebras was originated by
Moody, Rao and Yokonuma in [MRY] and [EM], where they constructed a vertex
operator representation in a homogeneous realization. The principal realization
was later given in [B1]. Both realizations were unified and substantially 
generalized in [L] and [BB].

As the first application of this representation theory, one may use the vertex operator
realizations to construct hierarchies of soliton equations as it was done in [B2],
[ISW] and [IT].

The modules for toroidal Lie algebras introduced in these papers have weight
decompositions with finite-dimensional weight spaces and are bounded, but
do not possess a unique highest weight. To explain this, we consider a $\Z$ grading
of $\g$ by degree with respect to $t_0$, which is declared to be a special variable.
The subalgebra of elements of degree $0$ in this $\Z$ grading, is very close to 
an $N$-toroidal Lie algebra. For this subalgebra we may consider an irreducible
module $T$, which is a ``toroidal'' module for  $(\R_0 \ot \dg) \oplus \K_0$
(a multi-variable analog of a loop module), and is a tensor module for $\D_0$.
We let the elements of positive degree act on $T$ trivially, and then induce $T$
to the module over the whole toroidal Lie algebra. This induced module has a unique
irreducible quotient that may be alternatively studied via the explicit vertex operator 
constructions. This approach to the representation theory of toroidal Lie algebras
was laid out in [BB]. As we see, instead of a one-dimensional highest weight space, the 
whole of $T$ will be the ``top'' of the resulting bounded module. The space $T$
is infinite-dimensional, but has a $\Z^N$ grading with finite-dimensional subspaces.

One serious problem with this representation theory 
that has not been previously resolved is that the vertex operator realizations were 
constructed not for the full toroidal algebra, but only
for its subalgebra $\g^* =(\R \ot \dg) \oplus \K \oplus \D^*$, 
where the derivations in $t_0$ are missing from the derivation part $\D$. 
Plausible candidates for representing the missing part yielded extremely messy cocycles 
with values in a certain complicated completion of $\K$ [MRY], [L]. 

This was leaving the whole picture incomplete and unsatisfactory from the physics
perspective, because the missing derivations are responsible for the energy-momentum
tensor for these modules. The main goal of the present paper is to resolve this 
problem and construct a class of representations for the full toroidal Lie algebras.

Often the representation theory of Lie algebras is used for the construction of 
the vertex operator algebras. In our case it is the opposite -- the representation
theory is developed using the machinery of the vertex operator algebras.
This has been done in [BBS] for the subalgebra $\g^*$ of the toroidal Lie algebra $\g$.  
The VOA that controls the representation theory of $\g^*$ is a tensor product of three
fairly well-known VOAs -- the affine $\dg$ VOA $V_{\wdg}$, the affine $\gl_N$ VOA
$V_{\wgl}$ and a sub-VOA $V_\hyp^+$ of a hyperbolic lattice VOA. 

After a careful analysis using the methods of [BB], it became clear that the irreducible
modules for $\g^*$ do not admit the action of the full toroidal Lie algebra $\g$, and
thus it is necessary to enlarge the representation space in order to get 
a module for $\g$. A natural guess is that this enlarged space should be again a VOA or
a VOA module. It turns out that the missing ingredient is a VOA that corresponds
to the twisted Heisenberg-Virasoro algebra $\HVir$, and the full toroidal VOA is 
a tensor product of four VOAs:
$$ V_\tor = V_\wdg \ot V_\hyp^+ \ot V_{\wsl} \ot V_\HVir $$
with certain conditions on the central charges of these VOAs.

The twisted Heisenberg-Virasoro Lie algebra has a Virasoro subalgebra and a Heisenberg
subalgebra, but the natural action of the Virasoro on the Heisenberg subalgebra is
twisted with a cocycle (see Section 2.4 for the precise definition). The representation
theory of $\HVir$ has been studied by Arbarello et al. in [ACKP]. However one
special case, namely when the central charge of the Heisenberg subalgebra is zero,
was not fully investigated in that paper. It happens that this is precisely the case
we need for the toroidal VOA. 
The structure of the irreducible modules for $\HVir$ with the trivial action of the 
center of the Heisenberg subalgebra has been determined in [B3].

Using these ingredients we can easily write down the characters of the toroidal VOA
and of its modules. This leads to the following open problem: while the explicit
expressions for the characters of irreducible modules are known, there is no Weyl-type
character formula for the toroidal Lie algebras. Obtaining such a formula may lead
to interesting number-theoretic identities. 

Toroidal Lie algebras are related to the class of extended affine Lie algebras.
These Lie algebras have been extensively studied during the last decade (see [AABGP]
and references therein). The main features of an extended affine Lie algebra
is that it is graded by a finite root system and possesses a non-degenerate
symmetric invariant bilinear form. The full toroidal Lie algebra $\g$ does not
possess a non-degenerate invariant form, but its subalgebra
$$ \gdiv = (\R \ot \dg) \oplus \K \oplus \D_\div$$
does. Here $\D_\div$ is the subalgebra of the divergence-free vector fields. 
Using the theory for the full toroidal algebra $\g$, we are able to construct
irreducible representations for its important subalgebra $\gdiv$ as well.

The vertex operator algebra that controls the representation theory of $\gdiv$ is
a tensor product of an affine VOA $V_\wdg$ at level $c$, a sub-VOA of a hyperbolic
lattice VOA $V_\hyp^+$, and a Virasoro VOA at level $\wcL$. The condition on the
central charges that we get here is
$$ 2(N+1) + {c \dim \dg \over c + h^\vee} + \wcL = 26,$$
which has a striking resemblance to the formula for the critical dimension
in the bosonic string theory.

Another interesting fact is that when $N = 12$, we get an exceptional module
for the Lie algebra $\D_\div \oplus \K$. Only for this value of $N$ we can
represent $\D_\div \oplus \K$ just on a hyperbolic lattice sub-VOA $V_\hyp^+$,
and the structure of the module becomes exceptionally simple. The character
of this module is given by the $-24$-th power of the Dedekind $\eta$-function
and has nice modular properties.

We should mention that the class of the modules for the full toroidal Lie algebra
constructed in this paper is not exhaustive, but could be described as a toroidal
counterpart of the level 1 representations for affine Lie algebras, so more
research remains to be done.

 The structure of the paper is the following. In Section 1 we give the definition
of the toroidal Lie algebras. In Section 2 we recall the definition and
the properties of VOAs and construct the vertex operator algebras corresponding
to the twisted Heisenberg-Virasoro algebra using the technique of the vertex Lie 
algebras. In Section 3 we describe
the tensor factors of the toroidal VOAs -- the affine VOA $V_\wdg$, a sub-VOA
of a hyperbolic lattice VOA $V_\hyp^+$, and the twisted $\wgl$-Virasoro VOA
$V_{\wgl-\Vir}$. We conclude Section 3 with the definition of the toroidal VOA.
In Section 4 we state and prove our main result -- every module for the toroidal VOA
is a module for the toroidal Lie algebra. In the last Section we describe the structure
of the irreducible modules for the full toroidal Lie algebra, as well as for its
subalgebra $\gdiv$. We conclude the paper with the construction of an exceptional
module for $\D_\div \oplus \K$ which is possible only when $N = 12$.

\

{\bf Acknowledgements:} I am grateful to Stephen Berman for the stimulating discussions
and encouragement.

\

{\bf 1. Toroidal Lie algebras.}

\

Toroidal Lie algebras are the natural multi-variable generalizations of
affine Lie algebras. 
In this review of the toroidal Lie algebras we follow the work [BB].
Let $\dg$ be a simple finite-dimensional Lie algebra
over $\C$ with a non-degenerate invariant bilinear form $( \cdot | \cdot )$
and let $N \geq 1$ be an integer. 
We consider the Lie algebra $\R \ot \dg$
of maps of an $N+1$ dimensional torus into $\dg$, where  
$\R = \C [t_0^\pm, t_1^\pm, \ldots, t_N^\pm]$ is the algebra of functions 
on a torus (in the Fourier basis). The universal central extension
of this Lie algebra may be described by means of the following construction which is
due to Kassel [Kas]. Let $\Omega_\R$ be the space of $1$-forms on a torus:
$\Omega_\R = \mathop\oplus\limits_{p = 0}^N \R dt_p$. We will choose the
forms $\{ k_p = t_p^{-1} dt_p | p = 0, \ldots, N \}$ as a basis of this
free $\R$ module. There is a natural map $d$ from the space of functions
$\R$ into $\Omega_\R$: $d(f) = \sum\limits_{p = 0}^N {\d f \over \d t_p} dt_p
= \sum\limits_{p = 0}^N t_p {\d f \over \d t_p} k_p$. The center $\K$ for the
universal central extension $(\R \ot \dg) \oplus \K$ of $\R \ot \dg$
is realized as
$$ \K = \Omega_\R / d(\R) , $$
and the Lie bracket is given by the formula
$$[f_1(t) g_1, f_2(t) g_2] = f_1(t) f_2(t) [g_1, g_2] + (g_1| g_2) f_2 d(f_1).$$
Here and in the rest of the paper we will denote elements of $\K$ by the
same symbols as elements of $\Omega_\R$, keeping in mind the canonical
projection $\Omega_\R \rightarrow \Omega_\R / d(\R)$.

Just as in affine case, we add to $(\R \ot \dg) \oplus \K$ the algebra $\D$ of outer derivations 
$$\D = \mathop\oplus\limits_{p=0}^N \R d_p,$$
where $ d_p = t_p {\d \over \d t_p}$.
We will denote the multi-indices by bold letters 
$\r = (r_0, r_1, \ldots, r_N)$, etc.,
and by $\t^\r$ the corresponding
monomials $t_0^{r_0} t_1^{r_1} \ldots t_N^{r_N}$.

The natural action of $\D$ on $\R \ot \dg$ 
$$[\t^\r d_a, \t^\m g] = m_a \t^{\r+\m} g \eqno{(\Ldg)}$$
uniquely extends to the action on the universal central extension
$(\R \ot \dg) \oplus \K$ by
$$[\t^\r d_a, \t^\m k_b] = m_a \t^{\r+\m} k_b + \delta_{ab}
\sum\limits_{p=0}^N r_p \t^{\r+\m} k_p . \eqno{(\Ldk)}$$
This corresponds to the Lie derivative action of the vector fields on 1-forms.

It turns out that there is still an extra degree of freedom in defining
the Lie algebra structure on $(\R \ot \dg) \oplus \K \oplus \D$.
The Lie bracket on $\D$ may be twisted with a $\K$-valued 2-cocycle:
$$[\t^\r d_a, \t^\m d_b] = m_a \t^{\r+\m} d_b - r_b \t^{\r+\m} d_a
+ \tau(\t^\r d_a, \t^\m d_b) . \eqno{(\Ldd)}$$
The space of these cocycles $H^2(\D,\K)$ is two-dimensional 
and is spanned by the following cocycles
$\tau_1$ and $\tau_2$:
$$\tau_1 (\t^\r d_a, \t^\m d_b) = m_a r_b \sum\limits_{p=0}^N
m_p \t^{\r+\m} k_p ,$$
$$\tau_2 (\t^\r d_a, \t^\m d_b) = r_a m_b \sum\limits_{p=0}^N
m_p \t^{\r+\m} k_p .$$
We will write $\tau = \mu \tau_1 + \nu \tau_2$. 
The resulting algebra (or rather a family of algebras) is called the
toroidal Lie algebra 
$$\g = \g(\mu,\nu) = (\R \ot \dg) \oplus \K \oplus \D.$$

Note that after adding the algebra of derivations $\D$, the center
of the toroidal Lie $\g$ becomes finite-dimensional with the basis 
$\{ k_0, k_1, \ldots, k_N \}$. This can be seen from the action
(\Ldk) of $\D$ on $\K$,  which is non-trivial.

In this paper we will consider only the multiples of the first cocycle
$\tau_1$, and we will be assuming $\nu = 0$ for most of our results here.

The toroidal Lie algebra $\g(\mu,\nu) = (\R \otimes \dg) \oplus \K \oplus \D$
has an important subalgebra $\gdiv (\mu)$ that has divergence free vector fields
as the derivation part:
$$ \gdiv(\mu) = (\R \otimes \dg) \oplus \K \oplus \Ddiv,$$
where
$$\Ddiv = \left\{ \sum_{p=0}^N f_p(\t) d_p \quad \bigg| \quad \sum_{p=0}^N t_p
{\d f_p \over \d t_p} = 0 \right\} .$$
The expression $i \sum\limits_{p=0}^N t_p  {\d f_p \over \d t_p}$ becomes
the divergence of a vector field in the angular coordinates $(x_0, \ldots, x_N)$
on a torus, where $t_j = e^{i x_j}$.

Note that the cocycle $\tau_2$ trivializes on $\Ddiv$, so we only get the
restriction of $\mu \tau_1$.

The importance of this subalgebra is explained by the fact that unlike the full
toroidal Lie algebra, $\gdiv$ is an extended affine Lie algebra [BGK],
i.e., $\gdiv$ has a non-degenerate symmetric invariant bilinear form. The
restrictions of this form to both $\R \otimes \dg$ and to $\Ddiv \oplus \K$ are
non-degenerate:
$$(\t^\r g_1 | \t^\m g_2 ) = \delta_{\r,-\m} (g_1 | g_2) , \quad g_1, g_2 \in \dg,$$
while the vector fields pair with the $1$-forms:
$$\big( \sum_{p=0}^N a_p \t^\r d_p | \t^\m k_q \big) = \delta_{\r,-\m} a_q . 
\eqno{(\pairDK)}$$
One can see that the above formula is ill-defined for the full $\D$,
since $d(\t^\m) = \sum\limits_{q=0}^N m_q \t^\m k_q$, being zero in $\K$, 
must be in the kernel
of the form. For the subalgebra $\Ddiv$ this is precisely the case since
$$ \big( \sum_{p=0}^N a_p \t^\r d_p |  \sum_{q=0}^N r_q \t^{-\r} k_q \big) = 
\sum\limits_{q=0}^N a_q r_q = 0.$$

All other values of the bilinear form are trivial:
$$(\R \otimes \dg | \Ddiv \oplus \K) = 0, \quad
(\Ddiv | \Ddiv) = 0, \quad (\K | \K) = 0 .$$
It is easy to verify that the resulting symmetric bilinear form is
invariant and non-degenerate.

The study representation theory of toroidal Lie algebras has begun in [MRY]
and [EM], with further developments in [B1], [L], [BB], [BBS]. In all of these
papers there was one common difficulty that has not been resolved --
the representations constructed there were not for the full toroidal
algebra $\g$, but only for a subalgebra
$$ \g^* = (\R \otimes \dg) \oplus \K \oplus 
\left( \mathop\oplus\limits_{p=1}^N \R d_p \right),$$
where the piece $\R d_0$ that corresponds to the toroidal energy-momentum 
tensor was missing. 
This left the theory in a somewhat incomplete form, and the goal of the present
paper is to construct a class of representations for the full
toroidal Lie algebra.  

\

\

{\bf 2. Vertex operator algebra associated with the twisted Heisenberg-
Virasoro Lie algebra.}

\

{\bf 2.1. Definitions and properties of a VOA.}

Let us recall the basic notions of the theory of the vertex operator algebras.
Here we are following [K2] and [Li].

{\bf Definition.} 
{\it
A vertex algebra is a vector space
$V$ with a distinguished vector
$\o$ (vacuum vector) in $V$, 
an operator $D$ (infinitesimal translation) on the space $V$, and a linear map $Y$ (state-field
correspondence)
$$\eqalign{
Y(\cdot,z): \quad V &\rightarrow (\End V)[[z,z^{-1}]], \cr
a &\mapsto Y(a,z) = \sum\limits_{n\in\Z} a_{(n)} z^{-n-1} 
\quad (\hbox{\it where \ } a_{(n)} \in \End V), \cr} $$
such that the following axioms hold:

\noindent
(V1) For any $a,b\in V, \quad a_{(n)} b = 0 $ for $n$ sufficiently large;

\noindent
(V2) $[D, Y(a,z)] = Y(D(a), z) = {d \over dz} Y(a,z)$ for any $a \in V$;

\noindent
(V3) $Y(\o,z) = \Id_V$;

\noindent
(V4) $Y(a,z) \o \in (\End V)[[z]]$ and $Y(a,z)\o |_{z=0} = a$ for any $a \in V$
\ (self-replication);

\noindent
(V5) For any $a, b \in V$, the fields $Y(a,z)$ and $Y(b,z)$ are mutually local, that is, 
$$ (z-w)^n \left[ Y(a,z), Y(b,w) \right] = 0, \quad \hbox{\it for \ } 
 n  \hbox{\it \ sufficiently large} .$$


A vertex algebra $V$ is called a vertex operator algebra (VOA) if, in addition, 
$V$ contains a vector $\omega$ (Virasoro element) such that

\noindent
(V6) The components $L(n) = \omega_{(n+1)}$ of the field
$$ Y(\omega,z) = \sum\limits_{n\in\Z} \omega_{(n)} z^{-n-1} 
= \sum\limits_{n\in\Z} L(n) z^{-n-2} $$
satisfy the Virasoro algebra relations:
$$  [ L(n) , L(m) ] = (n-m) L(n+m) + \delta_{n,-m} {n^3 - n \over 12} 
(\rank V) \Id, \quad \hbox{\it where \ } \rank V \in \C;  
\eqno{(\vir)}$$

\noindent
(V7) $D = L(-1)$;

\noindent
(V8) $V$ is graded by the eigenvalues of $L(0)$:
$V = \mathop\oplus\limits_{n\in \Z} V_n$ with $L(0) \big|_{V_n} = n \Id$.
}

This completes the definition of a VOA.

As a consequence of the axioms of the vertex algebra  
we have the  following important commutator formula:
$$\left[ Y(a,z_1), Y(b,z_2) \right] =
\sum_{n \geq 0} {1 \over n!}  Y(a_{(n)} b, z_2)
\left[ z_1^{-1} \left( {\d \over \d z_2} \right)^n
\delta \left( {z_2 \over z_1} \right) \right] . \eqno{(\comm)}$$
 As usual, the delta function is
$$ \delta(z) = \sum_{n\in\Z} z^n .$$
By (V1), the sum in the right hand side of the commutator formula
is actually finite.

All the vertex operator algebras 
that appear in this paper have the gradings by 
non-negative integers: $V = \mathop\oplus\limits_{n=0}^\infty V_n$.
In this case the sum in the right hand side of the commutator 
formula (\comm) runs from $n=0$ to $n = \deg(a) + \deg(b) -1$,
because 
$$\deg(a_{(n)} b) = \deg(a) + \deg(b) -n - 1 ,\eqno{(\degr)}$$ 
and the elements of negative degree vanish.

The commutator formula (\comm) may be written as the commutator
relations between the components of the vertex operators:
$$ [a_{(n)}, b_{(m)} ] =
\sum\limits_{j\geq 0} \pmatrix{ n \cr j \cr}
(a_{(j)} b)_{(n+m-j)} . \eqno{(\Borc)} $$
Equivalently,
$$ a_{(n)} b_{(m)} = b_{(m)} a_{(n)} +
\sum\limits_{j\geq 0} \pmatrix{ n \cr j \cr}
(a_{(j)} b)_{(n+m-j)} , \eqno{(\Bora)} $$
and also
$$ a_{(n)} b_{(m)} = b_{(m)} a_{(n)} -
\sum\limits_{j\geq 0} \pmatrix{ m \cr j \cr}
(b_{(j)} a)_{(n+m-j)} , \eqno{(\Borb)} $$

Another consequence of the axioms of a vertex algebra is the Borcherds'
identity:
$$\sum\limits_{j\geq 0} \pmatrix{m \cr j \cr}
(a_{(k+j)} b)_{(m+n-j)} c $$
$$= \sum\limits_{j\geq 0} 
(-1)^{k+j+1} \pmatrix{k \cr j \cr}
b_{(n+k-j)} a_{(m+j)} c +
\sum\limits_{j\geq 0} 
(-1)^j \pmatrix{k \cr j \cr}
a_{(m+k-j)} b_{(n+j)} c, \quad \quad  k,m,n \in \Z. \eqno{(\qas)} $$ 
We will particularly need its special case when $m=0$ and $k = -1$:
$$
(a_{(-1)} b)_{(n)} c = 
\sum\limits_{j\geq 0} b_{(n-j-1)} a_{(j)} c 
+ \sum\limits_{j\geq 0} 
a_{(-1-j)} b_{(n+j)} c, \quad \quad  k \in \Z. \eqno{(\qass)} $$

Let us list some other consequences of the axioms of a vertex algebra that we
will be using in the paper. It follows from V7 and V8 that 
$$ \omega_{(0)} a = D(a) \eqno{(\omd)} $$
and
$$ \omega_{(1)} a = \deg (a) a \quad \quad \hbox{\rm for \ } a
\hbox{\rm \  homogeneous}. \eqno{(\omdeg)} $$
The map $D$ is a derivation of the $n$-th product:
$$ D(a_{(n)} b) = (Da)_{(n)} b + a_{(n)} Db . \eqno{(\Dab)}$$
It could be easily derived from V2 that
$$\left( Da \right)_{(n)} = -n a_{(n-1)} \eqno{(\Da)}$$
and thus
$$ a_{(-1-k)} = {1\over k!} (D^k (a))_{(-1)} ,  \quad k \geq 0. 
\eqno{(\Dka)}$$

The last formula that we quote here is the skew-symmetry identity:
$$ a_{(n)} b = \sum_{j\geq 0} (-1)^{n+j+1} {1\over j!}
D^j (b_{(n+j)} a) . \eqno{(\skeww)}$$ 

\

{\bf 2.2. Tensor products of VOAs.}

The toroidal VOA that we introduce at the end of Section 3 is constructed
by taking a tensor product of three VOAs. Let us review here
the definition of the tensor product of two VOAs 
$\left( V^{\prime}, Y^{\prime}, \omega^{\prime}, \o \right)$
and $\left( V^{\prime\prime}, Y^{\prime\prime}, 
\omega^{\prime\prime}, \o \right)$
(the case
of an arbitrary number of factors is a trivial generalization).
The tensor product space $V = V^{\prime} \ot V^{\prime\prime}$ 
has the VOA structure under
$$ Y(a\ot b, z) = Y^{\prime} (a,z) \ot Y^{\prime\prime} (b,z), 
\eqno{(\Ytens)} ,$$
$$ \omega = \omega^{\prime} \ot \o 
+ \o \ot \omega^{\prime\prime},
\eqno{(\omtens)}$$
and $\o = \o \ot \o$ being the identity element.

It follows from (\omtens) that the rank of $V$ (see V6) is the sum of the ranks of the tensor factors. 

We will be later using the following simple lemma:

{\bf Lemma \ltens.} 
{\it Let $a,c \in V^{\prime}$, 
$b,d \in V^{\prime\prime}$. Then

\noindent
(i) $(a\ot \o)_{(-1)} (\o \ot b) = a \ot b.$     

\noindent 
(ii) $(a\ot \o)_{(n)} (\o \ot b) = 0$ for $ n \geq 0$.

\noindent
(iii) Suppose $a_{(j)} c = 0$ for $j\geq 0$. 
Then $(a \ot b)_{(n)} (c \ot d) = \sum\limits_{j\geq 0}
(a_{(-1-j)} c) \ot (b_{(n+j)} d) .$
}

{\it Proof.} Part (i) follows from V3 and V4. Part (ii) is 
a consequence of the commutativity of $Y(a\ot\o, z_1)$ and 
$Y(\o\ot b, z_2)$.
Part (iii) follows from (i), (ii) and (\qass).

\

{\bf 2.3. Vertex Lie algebras.}

 An important source of the vertex algebras is provided by the vertex Lie algebras.
In presenting this construction we will be following [DLM] (see also [P], [R],
[K2], [FKRW]).

Let $\L$ be a Lie algebra with the basis 
$\{ u(n), c(-1) \big| u\in\U, c\in\CC, n\in\Z \}$ ($\U$, $\CC$ are some index sets).
Define the corresponding fields in $\L [[z,z^{-1}]]$:
$$ u(z) = \sum_{n\in\Z} u(n) z^{-n-1}, \quad c(z) = c(-1) z^0, \quad 
u\in\U, c\in\CC .$$
Let $\F$ be a subspace in $\L [[z,z^{-1}]]$ spanned by all the fields
$u(z), c(z)$ and their derivatives of all orders.

{\bf Definition.}
{\it
 A Lie algebra $\L$ with the basis as above is called a vertex Lie algebra
if the following two conditions hold:

(1) for all $u_1, u_2 \in \U$,
$$ [u_1(z_1), u_2(z_2) ] = \sum\limits_{j=0}^n f_j(z_2)
\left[ z_1^{-1} \left( {\d \over \d z_2} \right)^j \delta \left(
{z_2 \over z_1} \right) \right], \eqno{(\vla)}$$
where $f_j(z) \in\F$ and $n$ depends on $u_1, u_2$,

(2) for all $c\in\CC$, the elements $c(-1)$ are central in $\L$.
}

\

This definition is a simplified version of the one from [DLM] and is 
not quite as general as the original definition, but it is sufficient 
for our purposes.

Let $\L^+$ be a subspace in $\L$ with the basis $\{ u(n) \big| u\in\U, n\geq 0 \}$
and let  $\L^-$ be a subspace with the basis 
$\{ u(n), c(-1) \big| u\in\U, c\in\CC, n<0 \}$. Then $\L = \L^+ \oplus \L^-$ and
$\L^+, \L^-$ are in fact subalgebras in $\L$.

The universal enveloping vertex algebra $V_\L$ of a vertex Lie algebra $\L$ 
is defined as an induced module
$$V_\L = \Ind_{\L^+}^\L (\C \o) = U(\L^-) \ot \o,$$
where $\C \o$ is a trivial 1-dimensional $\L^+$ module.

{\bf Theorem \voa.} ([DLM], Theorem 4.8) 
{\it
Let $\L$ be a vertex Lie algebra. Then

(a) $V_\L$ has a structure of a vertex algebra with the vacuum vector $\o$,
infinitesimal translation $D$ being a natural extension of 
the derivation of $\L$ given by 
$D(u(n))$  $=$ $-n u(n-1)$, $D(c(-1)) = 0$, $u\in\U$, $c\in\CC$,
and the state-field correspondence map $Y$ defined by the formula:
$$Y \left( a_1(-1-n_1) \ldots a_{k-1}(-1-n_{k-1}) a_k(-1-n_k) \o, z
\right) $$
$$ = 
:\left( {1\over n_1 !} \left( {\d \over \d z} \right)^{n_1} a_1 (z)
\right) 
\ldots
:\left( {1\over n_{k-1} !} \left( {\d \over \d z} \right)^{n_{k-1}} 
a_{k-1} (z) \right) 
\left( {1\over n_{k} !} \left( {\d \over \d z} \right)^{n_k} 
a_k (z) \right): \ldots : \quad 
, \eqno{(\Y)}$$
where $a_j \in \U, n_j \geq 0$ or $a_j \in\CC, n_j =0$.

(b) Any restricted $\L$ module is a vertex algebra module for $V_\L$.

(c) For an arbitrary character $\lambda: \CC \rightarrow \C$, the factor module
$$ V_\L (\lambda) =  U(\L^-) \o /  U(\L^-) \big< (c(-1) - \lambda(c)) 
\o \big>_{c\in\CC}$$
is a quotient vertex algebra.

(d) Any restricted $\L$ module in which $c(-1)$ act as $\lambda(c) \Id$, 
for all $c\in\CC$, is a vertex algebra module for $V_\L(\lambda)$.
}

In the formula (\Y) above,  
the normal
ordering of two fields $: a(z) b(z) :$ is defined as 
$$ : a(z) b(z) : = \sum\limits_{n < 0} a(n) z^{-n-1} b(z) + 
\sum\limits_{n \geq 0} b(z) a(n) z^{-n-1} .$$ 

\

{\bf 2.4. VOAs associated with the twisted Heisenberg-Virasoro algebra.}

 In the previous papers on the subject, no one has succeeded
in constructing a vertex operator representation
for the full toroidal Lie algebra.
However, Theorem 1.12 from [BB] predicts that such representations
should exist in this case as well. After a careful analysis using the 
methods
from [BB],
it became clear that irreducible $\g^*$ modules do not admit the action
of the full toroidal algebra $\g$, and their spaces should be enlarged
in order to make such an extension possible.
 The new  ingredient turned out to be the VOA associated with  
the twisted Heisenberg-Virasoro algebra which we describe
in this subsection.

 We define the twisted Heisenberg-Virasoro algebra $\HVir$ as a Lie 
algebra with the basis
$$\big\{ L(n), I(n), C_L, C_{LI}, C_I, \big| n \in \Z \big\}$$
and Lie bracket given by
$$ [ L(n) , L(m) ] = (n-m) L(n+m) + \delta_{n,-m} {n^3 - n \over 12} C_L 
,\eqno{(\LL)}$$
$$ [ L(n), I(m) ] = -m I(n+m) - \delta_{n,-m} (n^2 + n) C_{LI} ,
\eqno{(\LI)}$$
$$ [I(n), I(m) ] = n \delta_{n,-m} C_I, \eqno{(\II)}$$
$$[\HVir, C_L] = [\HVir, C_{LI}]  = [\HVir, C_I] = 0 .$$

This Lie algebra has a Heisenberg subalgebra and a Virasoro subalgebra
intertwined with the cocycle (\LI).
The twisted Heisenberg-Virasoro algebra $\HVir$ is the central extension 
of
the Lie algebra $\{ f(t) {d \over dt} + g(t) | f,g \in\C[t,t^{-1}] \}$ of
differential operators of order at most one. The corresponding
projection is given by $L(n) \mapsto -t^{n+1} {d \over dt}, $ \
$I(n) \mapsto t^n $.
The center of $\HVir$ is four-dimensional and is spanned by
$\{ I(0), C_L, C_{LI}, C_I \}$.

Irreducible highest weight representations for $\HVir$ have been studied 
by
Arbarello et al. [ACKP], where the structure of these modules 
is determined in case when the action of $C_I$ is non-zero. 
It turns out, however, that for our construction of the representations for 
the toroidal Lie algebras we need precisely the highest weight 
$\HVir$ modules in which $C_I$ acts as zero.
The irreducible modules of this type were studied in [B3], and we quote
the result from [B3] below.

We begin by recalling the standard construction of the Verma modules.

 Introduce a $\Z$ grading on $\HVir$ by deg$L(n) =$ deg$I(n) = n$ and
deg$C_L = \deg C_{LI} = \deg C_I = 0$,
and decompose $\HVir$ with respect to this grading:
$$\HVir = \HVir_- \oplus \HVir_0 \oplus \HVir_+.$$

Fix arbitrary complex numbers $h, h_I, c_L, c_{LI}, c_I$. Let $\C \o$ be a 
1-dimensional
$\HVir_0 \oplus \HVir_+$ module defined by $L(0) \o = h \o$, \
$I(0) \o = h_I \o$, \  $C_L \o = c_L \o$, \  $C_{LI} \o = c_{LI} \o$, \
$C_{I} \o = c_{I} \o$, \
$\HVir_+ \o = 0$.
As usual, the Verma module $M = M(h,h_I,c_L,c_{LI},c_I)$ is the induced 
module
$$ M(h,h_I,c_L,c_{LI},c_I) =
\hbox{\rm Ind}_{\HVir_0 \oplus\HVir_+}^{\HVir} (\C \o) \cong U(\HVir_-) 
\ot \o .$$

The module $M$ is $\Z$ graded by eigenvalues of the operator $L(0) - h$Id: 
\
$M = \mathop\oplus\limits_{n=0}^\infty M_n$ with $M_n = \{ v\in M | L(0) v 
= (n + h) v\}.$

 In order to understand the submodule structure of $M$, we need to study
singular vectors
in $M$.
A non-zero homogeneous vector $v$ in a highest weight $\HVir$ module is
called singular if $\HVir_+ v = 0$.
The Verma module $M(h,h_I,c_L,c_{LI},c_I)$ has a unique irreducible factor
which we denote $L(h,h_I,c_L,c_{LI},c_I)$.
\

{\bf Theorem \sing. ([B3], Theorem 1.)} 
{\it
Let $c_I = 0$ and $c_{LI} \neq 0.$

(a) If ${h_I \over c_{LI}} \notin \Z$ or  ${h_I \over c_{LI}} = 1$
then the $\HVir$ module $M(h,h_I,c_L,c_{LI},0)$ is irreducible.

(b) If ${h_I \over c_{LI}} \in \Z \backslash \{ 1 \}$ then
$M(h,h_I,c_L,c_{LI},0)$ possesses a singular vector $v\in M_n$,
where $n = \left| {h_I \over c_{LI}} - 1 \right|$.
The factor-module
$L(h,h_I,c_L,c_{LI},0) = M(h,h_I,c_L,c_{LI},0) / U(\HVir_-) v$ is 
irreducible
and its character is
$$ \char L(h,h_I,c_L,c_{LI},0) = (1 - q^n) \prod\limits_{k\geq 1} (1 - q^k)^{-2} .$$
}

\

 Using theorem \voa \ we can construct VOAs associated with the twisted 
Heisenberg-Virasoro algebra:

{\bf Theorem \hvmod .} 
{\it
Let $c_{LI} \neq 0$.

(a) The $\HVir$ module $L(0,0,c_L,c_{LI},c_I)$ is a simple vertex operator algebra.

(b) The $\HVir$ modules $M(h,h_I,c_L,c_{LI},0)$ and 
$L(h,h_I,c_L,c_{LI},0)$ are the VOA modules for $L(0,0,c_L,c_{LI},0)$.
}
 
{\it Proof.} First let us show that the twisted Heisenberg-Virasoro 
algebra is a vertex Lie algebra with $\U = \{ \om, I \}$ and
$\CC = \{ C_L, C_{LI}, C_I \}$, where we set
$\om(n) = L(n-1), C_L(-1) = C_L, C_{LI}(-1) = C_{LI}, C_I(-1) = C_I$.
Then the set 
$$\left\{ \om(n), I(n), C_L(-1), C_{LI}(-1), C_I(-1) \big| 
n\in\Z \right\}$$ 
is the basis of $\HVir$ compatible with the vertex structure.

Form the fields
$$\om(z) = \sum_{n\in\Z}\om(n) z^{n-1} = \sum_{n\in\Z} L(n) z^{n-2},$$  
$$I(z) = \sum_{n\in\Z} I(n) z^{n-1},$$
as well as the constant fields 
$C_L(z) = C_L, C_{LI}(z) = C_{LI}, C_I(z) = C_I$.  From the defining
relations (\LL),(\LI),(\II), we can derive the commutator relations
between $\om(z)$ and $I(z)$:
$$ \eqalign{
[\om(z_1), \om(z_2)] = \left( {\d \over \d z_2} \om(z_2) \right) \z 
&+ 2 \om(z_2) \zd \cr
&+ {C_L \over 12} \zddd ,} \eqno{(\LLz)}$$
$$ \eqalign{
[\om(z_1), I(z_2)] = \left( {\d \over \d z_2} I(z_2) \right) \z 
&+  I(z_2) \zd \cr
&- C_{LI} \zdd ,} \eqno{(\LIz)}$$
$$ [I(z_1), I(z_2)] = C_{I} \zd . \eqno{(\IIz)}$$
Observe that the expressions in the right hand sides are of the form
(\vla). Finally the elements $C_L, C_{LI}, C_I$ are central and we see
now that both properties of the vertex Lie algebra hold.

Let us prove that the irreducible $\HVir$ module $L(0,0,c_L,c_{LI},c_I)$
is a VOA. 
Let $\lambda$ be the character $\lambda(C_L) = c_L, 
\lambda(C_{LI}) = c_{LI}, \lambda(C_I) = c_I$. We are going to show
that $L(0,0,c_L,c_{LI},c_I)$ is a 
homomorphic image of the vertex algebra $V_\HVir (\lambda)$. 
As an $\HVir$ module, 
$V_\HVir (\lambda)$ is a factor of $U(\HVir)$ modulo the left ideal
generated by 
$$ \left\{ L(k), I(n), C_L - c_L\o, C_{LI} - c_{LI}\o, C_I - c_I \o
\big| k \geq -1, n \geq 0 \right\} .$$
The Verma module $M(0,0,c_L,c_{LI},c_I)$ on the other hand is a factor
of $U(\HVir)$ modulo the left ideal generated by 
$$ \left< L(n), I(n), C_L - c_L\o, C_{LI} - c_{LI}\o, C_I - c_I \o
\big| n \geq 0 \right> .$$  
Thus 
$$V_\HVir (c_L, c_{LI}, c_I) \cong M(0,0,c_L,c_{LI},c_I)/<L(-1)\o> 
\eqno{(\isom)}$$
(note that $L(-1)\o$ is a singular vector in $M(0,0,c_L,c_{LI},c_I)$).
Since $L(0,0,c_L,c_{LI},c_I)$ is a unique irreducible factor of
$M(0,0,c_L,c_{LI},c_I)$, we conclude that it is also a factor-module
of $V_\HVir (\lambda)$. Since the vertex operator algebra 
$V_\HVir (\lambda)$ is generated by the moments of the $\HVir$ fields
$\om(z)$ and $I(z)$, then every $\HVir$ submodule is a vertex algebra ideal
in $V_\HVir (\lambda)$. Thus every factor module of $V_\HVir(\lambda)$      
admits the structure of a (factor) vertex algebra, in particular the irreducible
module $L(0,0,c_L,c_{LI},c_I)$ becomes a simple vertex algebra. 
Note that $V_\HVir (\lambda)$ contains a Virasoro element $\om(-1)\o$,
so it is a VOA. The VOA structure
on $L(0,0,c_L,c_{LI},c_I)$ is still given by (\Y).

When $c_I = 0$, the vertex operator algebras 
$V_\HVir (c_L,c_{LI},0)$ and $L(0,0,c_L,c_{LI},0)$ are actually isomorphic
by Theorem \sing \ and (\isom). 
The second part of the theorem is now an
immediate consequence of Theorem \voa (d). This completes the proof
of the theorem.

\

Using the commutator formula (\comm) we derive from (\LIz) 
the following relations:

{\bf Lemma \omi.}

{\it
(i) $\omega_{(0)} I(-1) \o = D(I(-1)) \o, $

(ii) $\omega_{(1)} I(-1)\o = I(-1)\o, $

(iii) $\omega_{(2)} I(-1)\o = -2 c_{LI} \o, $

(iv) $\omega_{(j)} I(-1)\o = 0, $ for $j \geq 2$.
}

\

{\bf 3. Toroidal vertex operator algebra.}

\

In this section we discuss the tensor factors that will make up the
vertex operator algebra associated with the toroidal Lie algebra $\g$. 

For a subalgebra $\g^*$ of $\g$, the structure of the associated vertex
operator algebra and its modules was described in [BBS], however
until now the attempts to construct vertex operator representations for 
the full toroidal algebra $\g$ failed, which left the picture incomplete.

 Let us now describe three tensor factors of a toroidal VOA -- a 
sub-VOA of the hyperbolic lattice VOA, the affine VOA and the twisted 
$\wgl$-Virasoro VOA.  
The twisted $\wgl$-Virasoro VOA will be built using the twisted 
Heisenberg-Virasoro VOA introduced in the previous section.

\

{\bf 3.1. Hyperbolic lattice VOA.}

Here we present the construction of a hyperbolic lattice VOA. The general 
construction of a VOA corresponding to an arbitrary even lattice may be found 
in [FLM] or [K2]. 

Consider a hyperbolic lattice $\Hyp$, which is a free abelian group on $2N$
generators 
\break
$\{ u_i , v_i \}_{1 \leq i \leq N}$ with the symmetric bilinear
form
$$ ( \cdot | \cdot ) : \quad \Hyp \times \Hyp \rightarrow \Z ,$$
defined by 
$$ (u_i | v_j) = \delta_{ij} , \quad (u_i | u_j) = (v_i | v_j) = 0.$$
Note that the form $(\cdot | \cdot)$ is non-degenerate and $\Hyp$ is an 
even lattice, i.e., $(x | x) \in 2\Z$.

The construction of the VOA associated to $\Hyp$ proceeds as follows.

 First we complexify $\Hyp$:
$$ H = \Hyp \ot_{\Z} \C,$$
and extend $(\cdot | \cdot)$ by linearity on $H$. Next, we ``affinize'' $H$
by defining a Lie algebra
$\widehat H = \C[t, t^{-1}] \ot H \oplus \C K$
with the bracket
$$ [x(n), y(m)] = n (x | y) \delta_{n, -m} K, \quad x,y \in H, 
\quad [\widehat H, K] = 0.$$
Here and in what follows, we are using the notation $x(n) = t^n \ot x$. The algebra
$\widehat H$ has a triangular decomposition $\widehat H = \widehat H_-
\oplus \widehat H_0 \oplus \widehat H_+$, where 
$\widehat H_0 = < 1\ot H, K>$ and $\widehat H_\pm = 
t^{\pm 1} \C [t^\pm] \ot H$.

We also need a twisted group algebra of $\Hyp$, denoted by $\C[\Hyp]$, which
we now describe. The basis of $\C[\Hyp]$ is $\{ e^x | x \in \Hyp \}$, 
and the multiplication is twisted with the 2-cocycle $\epsilon$:
$$ e^x e^y = \epsilon(x,y) e^{x+y} , \quad x,y \in \Hyp , \eqno{(\twi)}$$
where $\epsilon$ is a multiplicatively bilinear map
$$\epsilon: \Hyp \times \Hyp \rightarrow \{ \pm 1 \},$$
defined on the generators by $\epsilon(v_i, u_j) = (-1)^{\delta_{ij}},
\epsilon(u_i, v_j) = \epsilon(u_i, u_j) = \epsilon(v_i, v_j) = 1,
\quad i,j = 1,\ldots, N$.

 We define the structure of $\widehat H_0 \oplus \widehat H_+$ module on 
$\C[\Hyp]$, letting $\widehat H_+$ act on $\C[\Hyp]$ trivially and
$\widehat H_0$ act by
$$ x(0) e^y = (x | y) e^y, \quad K e^y = e^y. \eqno{(\xoe)}$$

 Finally let $V_\hyp$ be the induced $\widehat H$ module:
$$ V_\hyp = \Ind_{\widehat H_0 \oplus \widehat H_+}^{\widehat H} \left(
\C[\Hyp] \right) .$$
This is the VOA attached to the lattice $\Hyp$. As a space $V_\hyp$ is 
isomorphic to the tensor product of the symmetric algebra $S(\widehat H_-)$
with the twisted group algebra $\C[\Hyp]$:
$$ V_\hyp = S(\widehat H_-) \ot \C[\Hyp] .$$

 The $Y$-map is defined on the basis elements of $\C[\Hyp]$ by
$$ Y(e^x, z) := \exp \left( \sum\limits_{j \geq 1} {x(-j) \over j} z^j \right)
  \exp \left( - \sum\limits_{j \geq 1} {x(j) \over j} z^{-j} \right)
e^x z^x , \eqno{(\Ye)}$$
where $e^x$ acts by twisted multiplication (\twi) and $z^x e^y = z^{(x|y)} e^y$.
For a general basis element $a = x_1 (-1-n_1) \ldots x_k(-1-n_k) \ot e^y$,
with $x_i, y \in \Hyp, n_i \geq 0$, one defines (cf. (\Y))
$$Y(a,z) = :\left( {1\over n_1 !} \left({\d \over \d z}
\right)^{n_1} x_1(z) \right) \ldots
\left( {1\over n_k !} \left({\d \over \d z}
\right)^{n_k} x_k(z) \right) Y(e^y,z) :, \eqno{(\Yw)}$$
where $x(z) = \sum\limits_{j\in\Z} x(j) z^{-j-1}$.
Note that $x_{(n)} = x(n)$ and sometimes the latter is more convenient
typographically.

The Virasoro element in $V_\hyp$ is $\omega = \sum\limits_{p=1}^N
u_p (-1) v_p(-1) \ot \o$, where $\o = e^0$ is the identity element
of $V_\hyp$. The rank of $V_\hyp$ is $2N$.

In the construction of the toroidal VOAs we would need not $V_\hyp$ itself,
but its sub-VOA $V_\hyp^+$:
$$  V_\hyp^+ = S(\widehat H_-) \ot \C[\Hyp^+] ,$$
where $\Hyp^+$ (resp. $\Hyp^-$) is the isotropic sublattice of $\Hyp$ 
generated by $\{ u_i, \quad 1\leq i \leq N \}$ (resp. 
$\{ v_i, \quad 1\leq i \leq N \}$). 
 One can verify immediately by inspecting
(\Ye) and (\Yw) that $V_\hyp^+$ is indeed a sub-VOA of $V_\hyp$. 
Also note that
the cocycle $\epsilon$ trivializes on $\C[\Hyp^+]$, making $\C[\Hyp^+]$
the usual (untwisted) group algebra.
The Virasoro element of $V_\hyp^+$ is the same as in $V_\hyp$, and so
the rank of $V_\hyp^+$ is also $2N$.

Let us describe a class of modules for $V_\hyp^+$. Consider the group
algebra $\C[H^+]$ of the vector space $H^+ = \Hyp^+ \ot_{\Z} \C.$
The space $S(\widehat H_-) \ot \C[H^+] \ot \C[\Hyp^-]$ has a 
structure of a VOA module for $V_\hyp^+$, where the action of 
$V_\hyp^+$ is still given by (\xoe),(\Ye) and (\Yw). 
Fix $\alpha \in \C^N, \beta\in\Z^N$. Then the subspace
$$ M_\hyp^+(\alpha,\beta) = S(\widehat H_-) \ot 
e^{\alpha \u + \beta \vv } \C[\Hyp^+]  \eqno{(\MHyp)}$$
in $S(\widehat H_-) \ot \C[H^+] \ot \C[\Hyp^-]$ is an irreducible
VOA module for $V_\hyp^+$. Here we are using the notations
$\alpha \u = \alpha_1 u_1 + \ldots +  \alpha_N u_N$,
$\m \u = m_1 u_1 + \ldots + m_N u_N$, etc.

The proof of the main result of this paper, Theorem \main \ will be comprised 
of some fairly deep calculations in a toroidal VOA.
The VOA $V_\hyp^+$ which is a tensor factor in a toroidal VOA will play
an essential role in these computations. We will present certain 
relations in $V_\hyp^+$ in a sequence of two simple lemmas. These 
relations will be used in the proof of Theorem \main.  

{\bf Lemma {\Hone}.}

{\it
(i) $\om_{(0)} \emu = D \emu = \sum\limits_{p=1}^N m_p u_p(-1) \emu, $

(ii) $\omega_{(j)} \emu = 0$ \quad for $j \geq 1$,

(iii)  $\omega_{(0)} u_a(-1) \emu = D (u_a(-1) \emu),$

(iv)  $\omega_{(1)} u_a(-1) \emu = u_a(-1) \emu,$

(v)  $\omega_{(j)} u_a(-1) \emu = 0$ \quad for $j \geq 2$,

(vi) $\omega_{(0)} v_a(-1) \emu = D (v_a(-1) \emu),$

(vii)  $\omega_{(1)} v_a(-1) \emu = v_a(-1) \emu,$

(viii)  $\omega_{(2)} v_a(-1) \emu = m_a \emu,$

(ix)  $\omega_{(j)} v_a(-1) \emu = 0$ \quad for $j \geq 3$.
}

{\it Proof.} Using (\omd) we get (iii), (vi) and the
first part of (i). The second part of (i) follows from
V2 and the fact that
$${d \over d z} Y(\emu,z) = \sum\limits_{p=1}^N m_p u_p(z) Y(\emu,z) 
= Y\left( \sum\limits_{p=1}^N m_p u_p(-1) \emu,z\right) .$$
The claims (iv), (vii) and (ii) with $j=1$ are the consequences of
 (\omdeg).
The statements (ii) with $j > 1$, (v) with $j > 2$ and (ix) follow
from (\degr). Finally, for (v) with $j=2$ and (viii) we recall
([FLM], (8.7.13)), that in a lattice VOA 
$$ [\omega_{(n+1)}, x(m)] = [L(n), x(m)] = -m x(m+n) .$$
Then for part (viii) we have
$$ \omega_{(2)} v_a (-1) \emu = [\omega_{(2)}, v_a(-1)]\emu +
v_a (-1) \omega_{(2)} \emu = v_a (0) \emu = m_a \emu ,$$  
while for part (v) with $j=2$ we get
$$ \omega_{(2)} u_a (-1) \emu = [\omega_{(2)}, u_a(-1)] \emu +
u_a (-1) \omega_{(2)} \emu = u_a (0) \emu = 0.$$ 
This completes the proof of the lemma.

\

{\bf Lemma \Htwo.} 

{\it
(i) $\eru_{(-1)} \emu = \ermu$,

(ii) $\eru_{(j)} \emu = 0$ \quad for $j \geq 0$,

(iii) $\eru_{(0)} v_a(-1) \emu = - r_a \ermu$,

(iv) $\eru_{(j)} v_a(-1) \emu = 0$ \quad for $j \geq 1$,

(v) $\eru_{(-1)} v_a(-1) \emu = v_a(-1) \ermu  - 
r_a (D\eru)_{(-1)}\emu$,

(vi) $\eru_{(j)} \omega_{(n)} =  \omega_{(n)} \eru_{(j)}
- (D \eru)_{(n+j)} $,

(vii)  $\left( u_s(-1) \eru \right)_{(0)} v_a(-1) \emu = 
\delta_{s,a} (D \eru)_{(-1)} \emu - r_a u_s(-1) \ermu$,

(viii) $\left( u_s(-1) \eru \right)_{(1)} v_a(-1) \emu = 
\delta_{s,a} \ermu$,

(ix) $\left( u_s(-1) \eru \right)_{(j)} v_a(-1) \emu = 0$
\quad for $j \geq 2$,

(x) $\eru_{(0)} D(v_a (-1)\emu) = -r_a D(\ermu)$, 

(xi) $\eru_{(1)} D(v_a (-1)\emu) = -r_a \ermu$. 
}

{\it Proof.}
Relation (i) follows from (\Ye), while (ii) is a consequence of the
commutativity of $Y(\eru,z_1)$ and $Y(\emu,z_2)$.

Using the commutator formula (\Borb) we get
$$ \eru_{(n)} v_a(-1) = v_a(-1) \eru_{(n)} -
\sum\limits_{j\geq 0} \pmatrix{ -1 \cr j \cr}
\left( v_a(j) \eru\right)_{(n-1-j)} $$
$$= v_a(-1) \eru_{(n)} - (v_a(0) \eru)_{(n-1)}
=  v_a(-1) \eru_{(n)} - r_a \eru_{(n-1)},$$
from which (iii),(iv) and (v) immediately follow.

The identity (vi) is obtained by applying (\Borb) and Lemma
\Hone (i),(ii).

To show equalities (vii)--(ix), we apply (\qass):
$$\left( u_s(-1) \eru \right)_{(n)} v_a(-1) \emu $$
$$ = \sum_{j\geq 0} \eru_{(n-1-j)} u_s(j) v_a(-1) \emu
+ \sum_{j\geq 0} u_s(-1-j) \eru_{(n+j)} v_a(-1) \emu .$$
The expression $u_s (j) v_a(-1) \emu$ is non-zero only for
$j=1$, in which case $u_s (1) v_a(-1) \emu = \delta_{s,a} \emu$,
while  by (iii) and (iv) $\eru_{(n+j)} v_a(-1) \emu$ is non-zero
only for $n=j=0$. Thus we get
$$\left( u_s(-1) \eru \right)_{(n)} v_a(-1) \emu =
\delta_{s,a} \eru_{(n-2)} \emu - r_a \delta_{n,0} u_s(-1) \ermu .$$
The claims (vii), (viii) and (ix) follow from this equality.   
 
For the last two statements of the lemma, we use (\Dab):
$$\eru_{(0)} D(v_a (-1)\emu) = D \left( \eru_{(0)}
v_a (-1)\emu \right) - (D \eru)_{(0)} v_a (-1)\emu
= - r_a D(\ermu) .$$
On the last step we used part (iii) of the lemma and (\Da). 
Equality (xi) is obtained in a similar way using (iv) and (iii):
$$\eru_{(1)} D(v_a (-1)\emu) = D \left( \eru_{(1)}
v_a (-1)\emu \right) - (D \eru)_{(1)} v_a (-1)\emu$$
$$ = \eru_{(0)} v_a(-1) \emu = - r_a \ermu .$$
This completes the proof of the lemma.

\

{\bf 3.2. Affine VOAs.}

 The theory of affine VOA is by now standard, so we outline it very briefly.
A detailed exposition may be found in [Li] or [K2].

 Let $\dg$ be a finite-dimensional simple Lie algebra over $\C$ with the 
symmetric invariant bilinear form normalized by the condition 
$(\theta | \theta) = 2$ for the longest root $\theta$ of $\dg$.

 The untwisted affine algebra associated with $\dg$ is 
$$ \widehat \dg = \C[t_0, t_0^{-1}] \ot \dg \oplus \C k_0,$$ 
where $k_0$ is a central element and
$$[g_1(n) , g_2(m)] = [g_1, g_2](n+m) + n \delta_{n,-m} (g_1 | g_2) k_0 .$$
As before, we write $g(n)$ for $t_0^n \ot g$. The affine Lie algebra
$\wdg$ may be identified with a subalgebra in the toroidal Lie algebra
$\g$.

It is easy to see that $\wdg$ is a vertex Lie algebra with the set $\U$
being a basis of $\dg$ and $\CC = \{ k_0 \}$.

 Let $c \in \C$ be an arbitrary constant and let $\C \o$ be the one-dimensional
module for $\C[t_0] \ot \dg \oplus \C k_0$ where $\C[t_0] \ot \dg$
acts trivially and $k_0$ acts as multiplication by $c$. 
By Theorem \voa \ the induced module
$$ V_\wdg (c) := \Ind_{\C[t_0] \ot \dg \oplus \C k_0}^\wdg (\C \o)
\cong U(t_0^{-1} \C[t_0^{-1}] \ot \dg) \ot \o $$
has a structure of a vertex algebra. 

If $c$ is not equal to the negative of the dual Coxeter number $h^\vee$ 
of $\dg$,
one can use the Sugawara construction to show that $V_\wdg (c)$ contains
a Virasoro element and thus turns into a VOA. The rank of this VOA is
$$ \rank V_\wdg (c) = {c \dim (\dg) \over c + h^\vee} .$$

Using (\omd), (\omdeg) and (\degr) we get the following

{\bf Lemma \omdd.}

{\it
(i) $\omega_{(0)} g(-1) \o = D(g(-1)) \o $,

(ii) $\omega_{(1)} g(-1) \o = g(-1) \o $,

(iii) $\omega_{(j)} g(-1) \o = 0$ for $j\geq 2$.
}

\

Finally, the irreducible quotient $L_{\wdg} (c)$ of $V_{\wdg} (c)$
is also a VOA of the same rank. By Theorem \voa (d) any highest weight
$\wdg$ module of level $c$ is a VOA module for $V_{\wdg} (c)$.

\

{\bf 3.3. $\wgl$ VOAs .}

 For the construction of the $(N+1)$-toroidal VOA we also need to consider the VOA associated with affine $\wgl$. 
The construction of this
VOA follows the general scheme of the previous subsection, with the
difference that $\glN$ is not simple, but only reductive.
The Lie algebra $\glN$ decomposes into a direct sum: $\glN = \slN \oplus \C I$,
where $I$ is an identity matrix. Let $\psi_1$ denote the projection on the traceless matrices and $\psi_2$ be the projection on the scalar matrices
in this decomposition:
$$\psi_2 (A) = {\tr (A) \over N} I, \quad \psi_1 (A) = A - \psi_2 (A), \quad
A \in \glN .$$
 Accordingly, we define the affine algebra $\wgl$
to be the direct sum of the affine algebra $\wsl = \C[t_0, t_0^{-1}] \ot
\slN \oplus \C C_1$ and a (degenerate) Heisenberg algebra $\Hei = 
\C[t_0, t_0^{-1}] \ot I \oplus \C C_I$. The Lie bracket in $\wgl$ is given
by
$$[g_1(n), g_2(m)] = [g_1,g_2](n+m) + n \delta_{n,-m} \left\{
tr(\psi_1(g_1)\psi_1(g_2)) C_1 + (\psi_2(g_1) | \psi_2(g_2)) C_I \right\} ,$$
where for the last term we will use normalization $(I | I) = 1$.

This can be rewritten as a commutator of the formal generating series:
$$\eqalign{
[g_1(z_1), g_2(z_2)] =& [g_1,g_2](z_2) \z \cr
&+ \left\{ tr(\psi_1(g_1)\psi_1(g_2)) C_1 
+  (\psi_2(g_1) | \psi_2(g_2)) C_I \right\}
\zd .  \cr} \eqno{(\ggz)}$$

 To construct a vertex algebra corresponding to $\wgl$, we will take a tensor
product of a vertex algebra for affine $\wsl$ and a Heisenberg vertex algebra. 
We outline the construction of the latter. 

The Heisenberg algebra $\C[t_0,t_0^{-1}] \ot I \oplus \C C_I$
with the Lie bracket given by (\II), is a vertex Lie algebra
with $\U = \{ I \}, \CC = \{ C_I \} .$
Let $c_I$ be an arbitrary constant in $\C$ and consider a one-dimensional
module $\C \o$ for the subalgebra $\C [t_0] \ot I \oplus \C C_I$ where 
$\C [t_0] \ot I$ acts on $\o$ trivially and $C_I$ acts as multiplication
by $c_I$. By theorem \voa, the induced module
$$V_\Hei (c_I) := \Ind_{\C [t_0] \ot I \oplus \C C_I}^{\Hei} (\C \o)
\cong U(t_0^{-1} \C [t_0^{-1}]\ot I) \ot \o ,$$
is a vertex algebra with
$$I(z) = \sum\limits_{j\in\Z} I(j) z^{-j-1}$$
and the map $Y$ defined by (\Y).

 The $\wgl$ vertex algebra, defined as a tensor product
$$ V_\wgl (c_1, c_I) := V_\wsl(c_1) \ot V_\Hei(c_I), $$
is a vertex operator algebra when
$c_1 \neq -N$ (the dual Coxeter number for $\slN$ is $N$) and $c_I \neq 0$.
In this case the irreducible quotient 
$L_\wgl(c_1,c_I) = L_\wsl (c_1) \ot V_\Hei(c_I)$ 
of $V_\wgl (c_1, c_I)$
is also a VOA. Any highest weight $\wgl$ module with $C_1, C_I$ acting as
multiplications by $c_1, c_I$, is a VOA module for $V_\wgl (c_1, c_I)$.

The following relations in $V_\wgl (c_1, c_I)$ are the 
consequences of (\ggz) and the commutator formula (\comm):

{\bf Lemma \EE.}

{\it
(i) $E_{ab} (-1)_{(0)} E_{cd} (-1) \o = \delta_{bc} E_{ad}(-1) \o - 
\delta_{ad} E_{cb}(-1) \o ,$

(ii) $E_{ab} (-1)_{(1)} E_{cd} (-1) \o = \delta_{ad} \delta_{bc} c_1 \o 
+ \delta_{ab} \delta_{cd} \left( {c_I \over N^2} - {c_1 \over N} 
\right) \o , $

(iii) $E_{ab} (-1)_{(j)} E_{cd} (-1) \o = 0, $ for $j\geq 2$.
}

\

{\bf 3.4. Twisted $\wgl$ - Virasoro VOAs.}

We define the twisted $\wgl$ - Virasoro VOA $V_{\wgl-\Vir}$ as the tensor product
of affine $\wsl$ VOA with the twisted Heisenberg-Virasoro VOA introduced in
Section 2.4:
$$ V_{\wgl-\Vir} (c_1, c_L, c_{LI}, c_I) := V_\wsl (c_1) \ot 
V_{\HVir} (c_L, c_{LI},c_I) .$$

The rank of this VOA is
$$ \rank V_{\wgl-\Vir} (c_1, c_L, c_{LI}, c_I) = 
{ c_1 (N^2 -1) \over c_1 + N } + c_L .$$

 We turn $V_{\wgl-\Vir}$ into a $\wgl$ module by identifying the symbol $I$
from the Heisenberg-Virasoro part with the identity matrix in $\glN$: 
$$ \sum\limits_{j\in \Z} A(j) z^{-j-1} \mapsto
Y_\wsl ( \psi_1 (A) (-1) \o, z) + Y_{\HVir} (\psi_2 (A) (-1) \o , z) =$$
$$ = \sum\limits_{j\in \Z} \psi_1 (A)(j) z^{-j-1} + {\tr(A) \over N}  
\sum\limits_{j\in \Z} I(j) z^{-j-1}, \quad {\rm for \ } A \in \glN .$$ 

 We may view the vertex algebra $V_\wgl (c_1, c_I)$ as a sub-vertex algebra
of the VOA 
\break
$V_{\wgl-\Vir} (c_1, c_L, c_{LI},c_I)$.

The VOA $V_{\wgl-\Vir} (c_1, c_L, c_{LI},c_I)$ has a simple factor VOA
$L_{\wgl-\Vir} (c_1, c_L, c_{LI},c_I) \cong  L_\wsl (c_1) \ot 
L_{\HVir} (0,0,c_L, c_{LI},c_I)$.

Let $L_\wsl (\lambda_1, c_1)$ be an irreducible $\wsl$ module of the
highest weight $(\lambda_1, c_1)$, where $\lambda_1$ is a linear functional
on the Cartan subalgebra of $\slN$. The tensor product 
$$ L_\wsl (\lambda_1, c_1) \ot L_{\HVir} (h,h_I,c_L, c_{LI},c_I)$$
is an irreducible VOA module for $V_{\wgl-\Vir}(c_1, c_L, c_{LI},c_I)$. 

We will later need the following relations in the twisted $\wgl$-Virasoro VOA:

{\bf Lemma \omE.}

{\it
(i) $\omega_{(0)} E_{ab} (-1) \o = D(E_{ab}(-1)) \o$,

(ii) $\omega_{(1)} E_{ab} (-1) \o = E_{ab}(-1) \o$,

(iii) $\omega_{(2)} E_{ab} (-1) \o = - \delta_{ab} 
{2 c_{LI} \over N} \o$,

(iv) $\omega_{(j)} E_{ab} (-1) \o = 0$ for $j \geq 3$.
}

{\it Proof.} Here (i) follows from (\omd), (ii) from (\omdeg),
and (iv) from (\degr). Let us establish (iii).
We project $E_{ab} (-1)\o$ into the affine $\wsl$ VOA and the 
twisted Heisenberg-Virasoro VOA: $E_{ab} (-1)\o = 
 \psi_1(E_{ab})(-1)\o + \delta_{ab} {1\over N}
I(-1) \o$. Applying Lemma \omdd (iii) and Lemma \omi (iii) we get  
$$\omega_{(2)} \psi_1(E_{ab}) (-1) \o = 0 ,$$
$$\omega_{(2)} I(-1) \o = -2 c_{LI} \o ,$$
and the claim (iii) follows.   

\

The main result of this paper is that for certain values of the 
parameters, the tensor product of VOAs
$$ V_{\tor} = V_\wdg (c) \ot V_\hyp^+ \ot V_{\wgl-\Vir} (0, c_L, c_{LI},0) \eqno{(\T)}$$
is a module for the toroidal Lie algebra $\g(\mu,0)$. The modules for this VOA
naturally inherit the structure of $\g(\mu,0)$ modules. 

 \

{\bf 4. Main theorem.}

\

The goal of the present paper is to show that the toroidal VOA (\T) has a structure  
of a module over the toroidal Lie algebra $\g (\mu, 0)$. 
In order to establish the correspondence between $V_{\tor}$ and $\g$, we consider the 
following fields in $\g [[z,z^{-1}]]$:
$$k_0 (\r,z) = \sum_{j\in\Z} t_0^j \t^\r k_0 z^{-j} ,$$
$$k_a (\r,z) = \sum_{j\in\Z} t_0^j \t^\r k_a z^{-j-1} , 
\quad a = 1,\ldots, N ,$$
$$g(\r,z) = \sum_{j\in\Z} t_0^j \t^\r g z^{-j-1} ,
\quad g \in \dg,$$
$$d_a (\r,z) = \sum_{j\in\Z} t_0^j \t^\r d_a z^{-j-1} ,
\quad a = 1,\ldots, N ,$$
$$\td (\r,z) = \sum_{j\in\Z} \left( - t_0^j \t^\r d_0 + 
\mu (j+{1\over 2}) t_0^j \t^\r k_0 \right) z^{-j-2} .$$
Here $\r\in\Z^N$ and $\t^\r = t_1^{r_1} \ldots t_N^{r_N}$.
In the last field we added the term $\mu (j+ {1\over 2})  t_0^j \t^\r k_0 $
in order to make the Lie bracket compatible with the VOA structure.
This is analogous to choosing the ``correct'' basis in a vertex Lie algebra.

The moments of the above fields span the toroidal algebra and the Lie bracket structure 
of $\g$ may be encoded in the commutator relations involving these fields. 
The commutators involving the first four fields only,
$k_0(\r,z), k_a(\r,z), g(\r,z), d_a(\r,z)$, are given in [BBS] and
will not be reproduced here. Let us write down only those commutators
that involve the last field $\td(\r,z)$. These relations are derived 
from (\Ldg),(\Ldk) and (\Ldd): 

$$
 [\td (\r, z_1), k_0(\m,z_2)] = \sum_{p=1}^N m_p k_p(\r+\m,z_2) 
\z , 
\eqno{(\dko)} $$
$$[\td (\r, z_1), k_a(\m,z_2)] = \dzb \left(  k_a(\r+\m,z_2) 
\z \right), \eqno{(\dka)}$$
$$ [\td (\r, z_1), g(\m,z_2)] = \dzb \left(  g(\r+\m,z_2) 
\z \right), \eqno{(\dgg)}$$
$$ [\td (\r, z_1), d_a(\m,z_2)] = \dzb \left(  d_a(\r+\m,z_2) 
\z \right)
- r_a \td(\r+\m,z_2) \z $$
$$ - \mu r_a \dzb \left( k_0(\r+\m,z_2) \zd \right)$$
$$ - \mu r_a \dzb \left( \sum_{p=1}^N r_p k_p(\r+\m,z_2) \z \right)
, \eqno{(\dda)} $$
$$ 
[\td (\r, z_1), \td(\m,z_2)] =  \td(\r+\m,z_2) \zd $$
$$ + \dzb \left(  \td(\r+\m,z_2) 
\z \right) $$
$$+ {\mu \over 2} \dzbb \left( k_0(\r+\m,z_2) \zd \right)$$
$$ + {\mu \over 2} \dzb \left( k_0(\r+\m,z_2) \zdd \right) $$
$$ + {\mu \over 2} \sum_{p=1}^N (r_p - m_p) \dzb \left( k_p(\r+\m,z_2) 
\zd \right)
. \eqno{(\ddo)}$$

Before we state the main theorem of the paper, we recall a result from [BBS], where the 
representation theory of a subalgebra $\g^*$ of the toroidal Lie 
algebra $\g$ is developed using the VOA approach.

{\bf Theorem \bbs.} ([BBS], Theorem 5.1.)
{\it
Let $M_{\wdg}$ be a module for affine vertex algebra $V_{\wdg}(c)$ with $c\neq 0$. 
Let $V_\hyp^+$ be a sub-VOA of 
a hyperbolic lattice VOA and let $M_{\wgl}$ be a module
for the vertex algebra $V_{\wgl} (c_1, c_I)$ . Then the tensor product
$$   M_{\wdg} \ot V_\hyp^+ \ot M_{\wgl} $$
is a module for the subalgebra $\g^*(\mu,\nu)$ of the toroidal Lie algebra
$\g(\mu,\nu)$, where the cocycle $\tau = \mu \tau_1 + \nu \tau_2$
has $\mu = {1 - c_1 \over c}$ and 
$\nu = {c_1 \over c N} - {c_I \over c N^2}$.
The action of the Lie algebra $\g^*(\mu,\nu)$ is given by the vertex 
operators:
$$k_0(\r,z) \mapsto c Y(\o\ot \eru \ot \o, z) , \eqno{(\ko)}$$
$$k_a(\r,z) \mapsto c Y(\o\ot u_a(-1)\eru \ot \o, z) , 
\eqno{(\ka)}$$
$$g(\r,z) \mapsto  Y(g(-1)\ot \eru \ot \o, z) , 
\eqno{(\ggg)}$$
$$d_a(\r,z) \mapsto  Y(\o\ot v_a(-1)\eru \ot \o, z) +
\sum_{p=1}^N r_p Y(\o\ot \eru \ot E_{pa}(-1),z). 
\eqno{(\da)}$$
}
  
In this paper we are taking $c_1 = c_I = 0$, so the corresponding
cocycle $\tau = \mu \tau_1 + \nu \tau_2$ has $\mu = {1\over c}$
and $\nu = 0$.

By a slight abuse of notations we will drop, from now on, the symbols $\o$ 
from the tensor products, e.g., write $\eru \ot E_{pa}(-1)$
instead of  $\o\ot \eru \ot E_{pa}(-1)$.

We are now able to state the main result of the paper. 

\

{\bf Theorem \main.} 
{\it
Let $c, c_L, c_{LI}$ be complex numbers such that
$c\neq 0, c \neq -h^\vee, c_{LI} = {N\over 2}$ and
$$ {c \dim \dg \over c + h^\vee} + 2N + c_L = 12. $$
Let also $c_1 = c_I = 0$.
Then the tensor product of the affine VOA, sub-VOA of a 
hyperbolic lattice VOA and the twisted $\wgl$-Virasoro VOA
$$V_\tor = V_\wdg(c) \ot V_\hyp^+ \ot V_{\wgl-\Vir}(0,c_L,c_{LI},0)$$
is a module for the toroidal Lie algebra $\g({1\over c},0)$.
The action of the fields
$k_0 (\r,z)$, $k_a(\r,z)$, $g(\r,z)$ and $d_a(\r,z)$ is given by
the formulas (\ko), (\ka), (\ggg), (\da) (however now in a different
representation space), whereas the field $\td(\r,z)$ is represented
by the vertex operator
$$\td(\r,z) \mapsto Y(\omega_{(-1)} \eru, z)
+ \sum_{p,s = 1}^N r_p 
Y( u_s(-1) \eru \ot E_{ps}(-1),z) ,
\eqno{(\tdd)}$$
where $\omega$ is the Virasoro element of the tensor product (the sum
of the Virasoro elements of $V_\wdg(c)$, $V_\hyp^+$ and
$V_{\wgl-\Vir}(0,c_L,c_{LI},0)$.
}

\

The following corollary can be immediately obtained from the theorem
by applying the principle of the ``preservation of identities''
([Li], Lemma 2.3.5).  

{\bf Corollary \tormod.} 
{\it
Let $c, c_L, c_{LI}$ satisfy the conditions of
Theorem \main. Let $M_\wdg$ be a VOA module for $V_\wdg(c)$, 
$M_\hyp^+$ be a VOA module for $V_\hyp^+$, and $M_{\wgl-\Vir}$ be
a VOA module for $V_{\wgl-\Vir}(0,c_L,c_{LI},0)$. Then 
$$ M_\wdg \ot M_\hyp^+ \ot M_{\wgl-\Vir}$$
is a module for the toroidal Lie algebra $\g({1\over c},0)$.
}

\
 
{\it Proof of Theorem \main.} 
We need to prove that the vertex operators in (\ko)--(\tdd)
satisfy the same commutator identities as the corresponding fields in the 
toroidal Lie algebra. Relations involving only the operators corresponding
to $k_0(\r,z), k_a(\r,z), g(\r,z)$ and $d_a(\r,z)$ follow from Theorem \bbs.
Our main tool for proving the remaining relations (\dko)--(\ddo) will be the 
commutator formula (\comm):
$$ \left[ Y(a,z_1), Y(b,z_2) \right] =
\sum\limits_{n=0}^{\deg(a) + \deg(b) - 1} {1\over n!} Y(a_n b, z_2)
\left[ z_1^{-1} \left({\d \over \d z_2} \right)^n \delta\left({z_1 \over z_2}
\right) \right] . \eqno{(\commm)}$$
Since our proof is somewhat lengthy, we organize it in four 
lemmas.

\

{\bf Lemma \Ydgg.} 
{\it
The vertex operators in (\tdd), (\ggg) 
representing $\td (\r,z)$ and $g(\m,z)$
satisfy the relation (\dgg).
}
 
{\it Proof.}
To verify that the vertex operators in (\tdd), (\ggg) satisfy (\dgg), 
we need to compute the $n$-th products
$$\left(
\om_{(-1)} \eru + \sum_{p,s =1}^N r_p  u_s (-1) 
\eru \ot E_{ps}(-1) \right)_{(n)} (g(-1) \ot \emu ) $$
for $n = 0, 1, 2$. The vertex operators 
$Y( u_s(-1) \eru \ot E_{ps}(-1),z_1)$ and
$Y(g(-1)\ot \emu , z_2)$  commute, and thus
$$\left(
\sum_{p,s =1}^N r_p  u_s (-1) 
\eru \ot E_{ps}(-1) \right)_{(n)} (g(-1) \ot \emu ) = 0,\quad {\rm for} \quad n\geq 0 .$$
Let us evaluate $\left(
\om_{(-1)}  \eru \right)_{(n)} (g(-1) \ot \emu )$ with the help 
of the Borcherds' identity (\qass):
$$\left(
\om_{(-1)}  \eru \right)_{(n)} (g(-1) \ot \emu ) $$
$$ = \sum_{j=0}^\infty \eru_{(n-j-1)} \om_{(j)}
(g(-1) \ot \emu) +
\sum_{j=0}^\infty  \om_{(-1-j)} \eru_{(j)} (g(-1) \ot \emu). \eqno
{(\Br)}$$  
Again because 
$Y( \eru ,z_1)$ and
$Y(g(-1)\ot \emu , z_2)$  commute, we get that 
\break
$\left( \eru  \right)_{(j)} 
\left(g(-1)\ot \emu \right) = 0$ for $j \geq 0$
and hence the last sum in $(\Br)$ equals $0$. For the first sum in 
(\Br) we write
$$ \om_{(j)} \left( g(-1) \ot \emu  \right) = 
(\om_{(j)} g(-1)) \ot \emu  + g(-1) \ot \om_{(j)} \emu  .
$$ 
By Lemma \Hone (i), (ii) we have $\om_{(0)} \emu = D \emu$, 
$\om_{(j)} \emu = 0$ for $j \geq 1$, 
and by Lemma \omdd,  $\om_{(0)} g(-1) = D g(-1)$,
$\om_{(1)} g(-1) = g(-1)$ and $\om_{(j)} g(-1) = 0$ for $j \geq 2$.
Thus (\Br) becomes
$$\left(\om_{(-1)}  \eru \right)_{(n)} 
(g(-1) \ot \emu ) = \eru_{(n-1)} D(g(-1) \ot \emu ) +
\eru_{(n-2)} (g(-1) \ot \emu ) . \eqno{(\omrumu)}
$$
Since $Y( \eru , z_1)$ commutes with 
$Y(g(-1) \ot \emu , z_2)$, we get that the first term in the
right hand side of (\omrumu) is non-zero for $n = 0$, while the second
term is non-zero for $n = 0, 1$.
Using Lemma \Htwo (i), we get that 
$$\left(\om_{(-1)} \eru \right)_{(1)} 
(g(-1) \ot \emu ) = g(-1) \ot \eru_{(-1)}\emu  =
g(-1) \ot \ermu ,  \eqno{(\omrumuone)}$$
$$\left(\om_{(-1)}  \eru  \right)_{(0)} 
(g(-1) \ot \emu ) =
\eru_{(-1)} D (g(-1) \ot \emu )
+ \eru_{(-2)}  (g(-1) \ot \emu ) = $$
$$ = \eru_{(-1)} D (g(-1) \ot \emu )
+ (D \eru)_{(-1)}  (g(-1) \ot \emu ) = 
D(g(-1) \ot \ermu ).  \eqno{(\omrumuo)}$$
  
Substituting (\omrumuone) and (\omrumuo) in the commutator formula    
(\commm) we get that 
$$\left[ Y(\om_{(-1)}\eru + \sum_{s,p=1}^N
r_s( u_p (-1) \eru \ot E_{ps} (-1)), z_1),
Y(g(-1)\ot\emu, z_2) \right] =$$
$$\left( {\d \over \d z_2} Y(g(-1)\ot\ermu, z_2) \right) \z
+ Y(g(-1)\ot\ermu, z_2) \zd ,$$
which corresponds to (\dgg). This completes the proof of the Lemma.

\

{\bf Lemma \Ydkk.} 
{\it
 The vertex operators in (\tdd), (\ko), (\ka) 
representing $\td (\r,z)$, \break
$k_0(\m,z)$ and $k_a (\m,z)$
satisfy the relations (\dko) and (\dka).
}

{\it Proof.}
We could perform the direct calculations for (\dko) and (\dka),
just as we did
in the previous lemma, however this computation may be
skipped altogether by 
referring to Lemma 2.1 of [BBS]. According to this lemma,
the relations (\dko) and (\dka) follow from (\dgg), which was established in Lemma \Ydgg.

\

{\bf Lemma \Ydda.}
{\it
 The vertex operators in (\tdd), (\da), (\ko), (\ka) 
representing $\td (\r,z)$, $d_a(\m,z)$, $k_0(\m,z)$ and $k_a (\m,z)$
satisfy the relation (\dda).
}

{\it Proof.} Again, to make use of the commutator formula (\commm), we need to compute the 
following $n$-th
products:
$$\left(
\om_{(-1)} \eru  + \sum_{p,s =1}^N r_p 
 u_s (-1) \eru \ot E_{ps}(-1) \right)_{(n)} 
\left( v_a(-1)\emu  +
\sum_{q=1}^N m_q  \emu \ot E_{qa}(-1) \right) \eqno{(\donda)}$$
for $n \geq 0$. 
We split (\donda) into four terms by expanding both sums. 
 In the four steps below, we handle each of these terms.

{\it Step 1.} Let us evaluate
$\left( \om_{(-1)}  \eru  \right)_{(n)} 
\left(  v_a(-1)\emu  \right). $
We apply the Borcherds' identity (\qass): 
$$\left( \om_{(-1)} \eru  \right)_{(n)} 
\left(  v_a(-1)\emu  \right) $$
$$ = \sum\limits_{j\geq 0} ( \eru )_{(n-1-j)}
\om_{(j)} (v_a(-1) \emu)  +
\sum\limits_{j\geq 0} \om_{(-1-j)}
\eru_{(n+j)}v_a(-1) \emu  .$$

Since the left factor in the $n$-th product (\donda) is of degree 2,
and the right factor is of degree 1, then by (\degr), we need only
to consider values $n = 0, 1, 2$.
We apply Lemma \Hone (vi)--(ix) to the first sum in the right hand side, and 
Lemma \Htwo (iii),(iv) to the second sum.

Let $n = 2$. In this case
$$\left( \om_{(-1)} \eru \right)_{(2)} 
\left(  v_a(-1)\emu  \right) 
= \eru_{(1)} D(v_a (-1)\emu) +
\eru _{(0)} v_a (-1)\emu +
\eru_{(-1)} ( m_a \emu) $$
$$= -r_a \ermu - r_a\ermu + m_a \ermu = (m_a - 2 r_a) \ermu .$$
In the above calculation we also used Lemma \Htwo (iii), (xi).

Next let $n=1$.  
$$\left( \om_{(-1)} ( \eru ) \right)_{(1)} 
\left(  v_a(-1)\emu  \right) 
=  \eru _{(0)} D(v_a (-1)\emu) +
\eru_{(-1)} v_a (-1)\emu +
\eru_{(-2)} ( m_a \emu) $$
$$ = - r_a  D(\ermu)  +  v_a(-1) \ermu 
- r_a  (D\eru)_{(-1)} \emu  +
m_a  (D\eru)_{(-1)} \emu  .$$
Here we used Lemma \Htwo (x), (v) and (\Dka).

Finally let $n=0$.
$$\left( \om_{(-1)} ( \eru ) \right)_{(0)} 
\left( v_a(-1)\emu  \right) $$
$$= \eru_{(-1)} D(v_a (-1)\emu) +
\eru_{(-2)} v_a (-1)\emu +
\eru_{(-3)} ( m_a \emu) 
+ \om_{(-1)} \eru_{(0)} v_a (-1)\emu$$
$$=  D(\eru_{(-1)} v_a (-1)\emu) 
+ {m_a \over 2}  (D^2\eru)_{(-1)}\emu 
- r_a \om_{(-1)} \ermu $$
$$=  D(v_a (-1)\ermu) 
- r_a  D((D \eru)_{(-1)} \emu) 
+ {m_a \over 2}  (D^2\eru)_{(-1)}\emu 
- r_a \om_{(-1)} \ermu . $$

To obtain the last two equalities we applied (\Dab), (\Dka) and
Lemma \Htwo (iii), (v).

The computations for Step 1 are now complete and we pass to

{\it Step 2.} We are going to evaluate 
$$\left( \om_{(-1)} \eru \right)_{(n)} 
\sum_{q=1}^N m_q \emu \ot E_{qa}(-1) , 
$$
for $n = 0,1,2$. Just as in Step 1, we use the Borcherds' identity
(\qass):
$$\left( \om_{(-1)} \eru \right)_{(n)} 
\sum_{q=1}^N m_q   \emu \ot E_{qa}(-1) $$
$$ = \sum\limits_{j\geq 0} \eru_{(n-1-j)}
\sum_{q=1}^N m_q \left(  (\om_{(j)} \emu) \ot E_{qa}(-1) 
+  \emu \ot \om_{(j)} E_{qa}(-1) \right)$$
$$+ \sum\limits_{j\geq 0} \om_{(-1-j)} 
\eru_{(n+j)} \sum_{q=1}^N m_q \emu\ot E_{qa} (-1) . $$
The second sum is equal to zero, while the first one yields (we use 
Lemma \Hone (i), (ii), Lemma {\omE} and the condition 
$c_{LI} = {N\over 2})$:
$$\left( \om_{(-1)} \eru \right)_{(n)} 
\sum_{q=1}^N m_q  \emu \ot E_{qa}(-1) $$
$$ = \eru_{(n-1)} \sum_{q=1}^N m_q  D(\emu) \ot 
E_{qa}(-1) +
\eru_{(n-1)} \sum_{q=1}^N m_q  \emu \ot 
D E_{qa}(-1)$$ 
$$+ \eru_{(n-2)} \sum_{q=1}^N m_q  \emu \ot E_{qa}(-1) 
- m_a {2 c_{LI} \over N} \eru_{(n-3)}\emu$$
$$ =  \eru_{(n-1)} \sum_{q=1}^N m_q D( \emu \ot 
E_{qa}(-1)) +
\eru_{(n-2)} \sum_{q=1}^N m_q  \emu \ot E_{qa}(-1)
- m_a  \eru_{(n-3)}\emu  . \eqno{(\Ib)}$$ 

Let $n=2$. Then the only non-zero term in the above expression
is the last term:
$$\left( \om_{(-1)} \eru \right)_{(2)} 
\sum_{q=1}^N m_q  \emu \ot E_{qa}(-1) 
= - m_a \ermu .$$

When $n=1$, the two last terms in (\Ib) are non-zero:
$$\left( \om_{(-1)} \eru \right)_{(1)} 
\sum_{q=1}^N m_q \emu \ot E_{qa}(-1) = \sum_{q=1}^N m_q \ermu\ot E_{qa}(-1) 
- m_a  (D\eru)_{(-1)}\emu .$$

For $n=0$ we get:
$$\left( \om_{(-1)} \eru \right)_{(0)} 
\sum_{q=1}^N m_q  \emu \ot E_{qa}(-1) $$
$$= \eru_{(-1)} \sum_{q=1}^N m_q D(\emu\ot E_{qa}(-1))
+ (D\eru)_{(-1)} \sum_{q=1}^N m_q \emu\ot E_{qa}(-1) 
- {m_a \over 2}  (D^2 \eru)_{(-1)}\emu $$
$$= \sum_{q=1}^N m_q D(\ermu\ot E_{qa}(-1))
- {m_a \over 2}  (D^2 \eru)_{(-1)}\emu .$$

{\it Step 3.}
The next term from (\donda) that we need to handle is
$$\left(
\sum_{p,s =1}^N r_p 
 u_s (-1) \eru \ot E_{ps}(-1) \right)_{(n)} 
v_a(-1)\emu .$$
Using Lemma \ltens (iii) we get
$$\left(
\sum_{p,s =1}^N r_p 
 u_s (-1) \eru \ot E_{ps}(-1) \right)_{(n)} 
v_a(-1)\emu $$
$$ = \sum\limits_{j\geq 0} \sum_{p,s =1}^N r_p 
 ((u_s (-1) \eru)_{(n+j)}  v_a(-1)\emu) \ot E_{ps}(-1-j).$$
Next we apply Lemma \Htwo (ix) and we see that the
expression above turns into zero for $n \geq 2$. Let us consider
now cases $n = 0, 1$. 

Let $n = 1$. Then using Lemma \Htwo (viii), (ix) we obtain
$$\left(
\sum_{p,s =1}^N r_p 
 u_s (-1) \eru \ot E_{ps}(-1) \right)_{(1)} v_a(-1)\emu 
 = \sum_{p=1}^N r_p \ermu \ot E_{pa}(-1) .$$

For $n = 0$ we have
$$\left(
\sum_{p,s =1}^N r_p 
 u_s (-1) \eru \ot E_{ps}(-1) \right)_{(0)} 
v_a(-1)\emu 
= \sum_{p=1}^N r_p   ((D \eru)_{(-1)} \emu) \ot E_{pa}(-1)$$
$$- r_a \sum_{p,s =1}^N r_p   u_s (-1) \ermu \ot E_{ps}(-1)
+ \sum_{p=1}^N r_p   \ermu \ot E_{pa}(-2) .$$

{\it Step 4.} 
We evaluate 
$$\left( \sum_{p,s =1}^N r_p 
 u_s (-1) \eru \ot E_{ps}(-1) \right)_{(n)} 
\left( \sum_{q=1}^N m_q  \emu \ot E_{qa}(-1) \right)$$
using Lemma \ltens (iii).
 Since we assume $c_1 = c_I = 0$ then by Lemma 
{\EE} we have 
\break
$E_{ps}(-1)_{(j)} E_{qa}(-1)$ $= 0$ for $j \geq 1$.
Thus the above expression is non-zero only for $n=0$, 
in which case we get:
$$\left( \sum_{p,s =1}^N r_p 
 u_s (-1) \eru \ot E_{ps}(-1) \right)_{(0)} 
\left( \sum_{q=1}^N m_q  \emu \ot E_{qa}(-1) \right)$$
$$ = \sum_{p,s =1}^N \sum_{q=1}^N r_p m_q  
(u_s(-1) \eru)_{(-1)} \emu
\ot E_{ps}(-1)_{(0)} E_{qa}(-1) $$
$$ = \sum_{p,s =1}^N  r_p m_s  u_s(-1) \ermu \ot E_{pa}(-1)
- r_a \sum_{q,s =1}^N m_q  u_s(-1) \ermu \ot E_{qs}(-1) $$  
$$ = \sum_{p=1}^N  r_p  (\eru_{(-1)} D\emu) \ot E_{pa}(-1)
- r_a \sum_{q,s =1}^N m_q  u_s(-1) \ermu \ot E_{qs}(-1) .$$ 

Now we combine the results of the four steps:
$$\left(
\om_{(-1)} \eru  + \sum_{p,s =1}^N r_p 
 u_s (-1) \eru \ot E_{ps}(-1) \right)_{(2)} 
\left(  v_a(-1)\emu  +
\sum_{q=1}^N m_q  \emu \ot E_{qa}(-1) \right)$$
$$ = (m_a - 2 r_a) \ermu - m_a \ermu  = - 2 r_a \ermu .$$
This corresponds to the term in (\dda) with the second derivative
of the delta function.

 Now let $n = 1$:
$$\left(
\om_{(-1)}  \eru  + \sum_{p,s =1}^N r_p 
 u_s (-1) \eru \ot E_{ps}(-1) \right)_{(1)} 
\left(  v_a(-1)\emu  +
\sum_{q=1}^N m_q  \emu \ot E_{qa}(-1) \right)$$
$$ =  v_a(-1) \ermu 
- r_a  D(\ermu)  - r_a  (D\eru)_{(-1)} \emu  +
m_a  (D\eru)_{(-1)} \emu  $$
$$+ \sum_{q=1}^N m_q \ermu\ot E_{qa}(-1) 
- m_a  (D \eru)_{(-1)}\emu 
+ \sum_{p=1}^N r_p  \ermu \ot E_{pa}(-1) $$
$$ =   v_a(-1) \ermu 
+ \sum_{p=1}^N (r_p + m_p)  \ermu \ot E_{pa}(-1)
- r_a  D(\ermu)  - r_a  (D\eru)_{(-1)} \emu  .$$
This matches the terms in (\dda) with the first derivative of the
delta function.

 Finally let $n=0$:
$$\left(
\om_{(-1)}  \eru  + \sum_{p,s =1}^N r_p 
 u_s (-1) \eru \ot E_{ps}(-1) \right)_{(0)} 
\left(  v_a(-1)\emu  +
\sum_{q=1}^N m_q  \emu \ot E_{qa}(-1) \right)$$
$$=  D(v_a (-1)\ermu) 
- r_a  D((D \eru)_{(-1)} \emu) 
+ {m_a \over 2}  (D^2\eru)_{(-1)}\emu 
- r_a \om_{(-1)} \ermu $$
$$+ \sum_{q=1}^N m_q D(\ermu\ot E_{qa}(-1))
- {m_a \over 2}  (D^2 \eru)_{(-1)}\emu 
+ \sum_{p=1}^N r_p   ((D \eru)_{(-1)} \emu) \ot E_{pa}(-1)$$
$$- r_a \sum_{p,s =1}^N r_p   u_s (-1) \ermu \ot E_{ps}(-1)
+ \sum_{p=1}^N r_p   \ermu \ot D(E_{pa}(-1)) $$
$$ + \sum_{p=1}^N  r_p  (\eru_{(-1)} D\emu) \ot E_{pa}(-1)
- r_a \sum_{q,s =1}^N m_q  u_s(-1) \ermu \ot E_{qs}(-1) $$  
$$ = D \left(
 v_a (-1)\ermu  + 
\sum_{p=1}^N (r_p + m_p)   \ermu \ot E_{pa}(-1) \right)$$
$$ - r_a \left(  \om_{(-1)}  \ermu 
+ \sum_{p,s =1}^N (r_p + m_p)   u_s (-1) \ermu \ot E_{ps}(-1) \right)
 - r_a  D((D \eru)_{(-1)} \emu)  .$$
This expression matches the terms with the delta function in (\dda)
and the Lemma is proved.

\

{\bf Lemma \Yddo.}
{\it
 The vertex operators in (\tdd), (\ko), (\ka) 
representing $\td (\r,z), k_0(\m,z)$ and $k_a (\m,z)$
satisfy the relation (\ddo).
}

{\it Proof.} We will verify (\ddo) using the commutator
formula (\commm). For this
we need to evaluate the following $n$-th products:
$$\left(
\om_{(-1)} \eru  + \sum_{p,s =1}^N r_p 
 u_s (-1) \eru \ot E_{ps}(-1) \right)_{(n)} 
\left(
\om_{(-1)} \emu  + \sum_{q,k =1}^N m_q 
 u_k (-1) \emu \ot E_{qk}(-1) \right). \eqno{(\dondo)}$$
We expand the sums in both factors and get four terms, which will be
dealt with in the four steps below.
Since both factors are of degree 2, we need to consider values
$n = 0,1,2,3$ (see (\commm)).

{\it Step 1.} The first term to be simplified is
$\left( \om_{(-1)}  \eru  \right)_{(n)} 
\left( \om_{(-1)}  \emu   \right)$.
We apply the Borcherds' identity (\qass):
$$\left(
\om_{(-1)}  \eru  \right)_{(n)} 
\left(
\om_{(-1)} \emu   \right)$$
$$ = \sum_{j\geq 0} \eru_{(n-j-1)} 
\om_{(j)} \om_{(-1)} \emu
+  \sum_{j\geq 0} \om_{(-1-j)} 
\eru_{(n+j)} \om_{(-1)} \emu. \eqno{(\ddos)}$$
For the first sum we use the relations in the Virasoro
algebra: 
$$\om_{(j)} \om_{(-1)} = \om_{(-1)} \om_{(j)}
+ (j+1) \om_{(j-2)} + 6 \delta_{j,3} \Id$$ 
(note that by
the assumption of the theorem, the rank of the VOA $V_\tor$ is 
equal to 12). We transform the second sum in (\ddos) using 
Lemma \Htwo (vi): 
$$\sum_{j\geq 0} \om_{(-1-j)} \eru_{(n+j)} \om_{(-1)} \emu$$
$$= \sum_{j\geq 0} \om_{(-1-j)} \om_{(-1)} \eru_{(n+j)} \emu
- \sum_{j\geq 0} \om_{(-1-j)} (D \eru)_{(n+j-1)} \emu
= - \delta_{n,0} \om_{(-1)} (D \eru)_{(-1)} \emu .$$
Thus (\ddos) becomes
$$\left(
\om_{(-1)}  \eru  \right)_{(n)} 
\left( \om_{(-1)}  \emu   \right)$$
$$ = \sum_{j\geq 0} \eru_{(n-j-1)}\om_{(-1)} \om_{(j)} 
\emu 
+  \sum_{j\geq 0} (j+1)\eru_{(n-j-1)}\om_{(j-2)}  \emu $$
$$+ 6 \eru_{(n-4)} \emu 
- \delta_{n,0} \om_{(-1)} (D\eru)_{(-1)} \emu .$$
This can be simplified further with the help of Lemma \Hone (i),(ii)
and Lemma \Htwo (vi):
$$\left(
\om_{(-1)}  \eru  \right)_{(n)} 
\left( \om_{(-1)}  \emu   \right)$$
$$ = \eru_{(n-1)} \om_{(-1)} D\emu 
+ \eru_{(n-1)} \om_{(-2)} \emu
+ 2 \eru_{(n-2)} \om_{(-1)} \emu
+ 3 (\eru)_{(n-3)} \om_{(0)} \emu$$
$$+ 6  \eru_{(n-4)} \emu 
- \delta_{n,0} \om_{(-1)} ( (D\eru)_{(-1)} \emu ) $$
$$ = \om_{(-1)} \eru_{(n-1)} D\emu - (D\eru)_{(n-2)} D\emu
+ \om_{(-2)} \eru_{(n-1)} \emu - (D\eru)_{(n-3)} \emu$$
$$+ 2 \om_{(-1)} \eru_{(n-2)} \emu - 2 (D\eru)_{(n-3)} \emu
+ 3 \eru_{(n-3)} D\emu$$
$$ + 6 \eru_{(n-4)} \emu
- \delta_{n,0} \om_{(-1)} ( (D\eru)_{(-1)} \emu ) $$
$$ = \om_{(-1)} \eru_{(n-1)} D\emu
+ \om_{(-2)} \eru_{(n-1)} \emu + 2 \om_{(-1)} \eru_{(n-2)} \emu
- \delta_{n,0} \om_{(-1)} ( (D\eru)_{(-1)} \emu ) $$
$$- (D\eru)_{(n-2)} D\emu  - 3 (D\eru)_{(n-3)} \emu
+ 3 \eru_{(n-3)} D\emu + 6 \eru_{(n-4)} \emu .$$
Now we consider particular values of $n$. Let $n=3$. In this case
the previous expression will simplify to just one term:
$$\left( \om_{(-1)} ( \eru ) \right)_{(3)} 
\left( \om_{(-1)} ( \emu  ) \right)
= 6 \ermu .$$
When $n=2$ we obtain using (\Dab) and (\Dka):
$$\left( \om_{(-1)} ( \eru ) \right)_{(2)} 
\left( \om_{(-1)} ( \emu  ) \right)
= - 3 (D\eru)_{(-1)} \emu
+ 3 \eru_{(-1)} D\emu + 6 \eru_{(-2)} \emu$$
$$ = 3 D \ermu .$$ 
For $n=1$ we get:
$$\left( \om_{(-1)} \eru \right)_{(1)} 
\left( \om_{(-1)} \emu \right)$$
$$= 2 \om_{(-1)} \ermu - (D\eru)_{(-1)} D\emu
- 3 (D\eru)_{(-2)} \emu
+ 3 \eru_{(-2)} D\emu + 6 \eru_{(-3)} \emu$$
$$ = 2 \om_{(-1)} \ermu + 2 (D\eru)_{(-1)} D\emu .$$
Finally, for $n=0$ we have
$$\left( \om_{(-1)} \eru \right)_{(0)} 
\left( \om_{(-1)} \emu \right)$$
$$= \om_{(-1)} \eru_{(-1)} D\emu
+ \om_{(-2)} \eru_{(-1)} \emu + 2 \om_{(-1)} \eru_{(-2)} \emu
- \om_{(-1)} (D\eru)_{(-1)} \emu $$
$$ - (D\eru)_{(-2)} D\emu  - 3 (D\eru)_{(-3)} \emu
+ 3 \eru_{(-3)} D\emu + 6 \eru_{(-4)} \emu $$
$$ = D(\om_{(-1)} \ermu) + {1\over 2} (D^2 \eru)_{(-1)} D\emu 
- {1\over 2} (D^3 \eru)_{(-1)} \emu .$$

{\it Step 2.} 
As in the previous step, to compute the $n$-th product
$$\left(
\om_{(-1)} \eru \right)_{(n)} 
\sum_{q,k =1}^N m_q  u_k (-1) \emu \ot E_{qk}(-1) $$
we use the Borcherds' identity (\qass), as well as (\omd), (\omdeg),
Lemma \Hone (v) and Lemma {\omE} (with $c_{LI}={N\over 2}$):
$$\left(
\om_{(-1)} \eru \right)_{(n)} 
\sum_{q,k =1}^N m_q  u_k (-1) \emu \ot E_{qk}(-1) $$
$$ = \sum_{j\geq 0} \eru_{(n-j-1)} 
\sum_{q,k =1}^N m_q  
\om_{(j)} (u_k (-1) \emu \ot E_{qk}(-1))$$
$$ = \eru_{(n-1)} \sum_{q,k =1}^N m_q  
\om_{(0)} (u_k (-1) \emu \ot E_{qk}(-1)) 
+ \eru_{(n-2)} \sum_{q,k =1}^N m_q  
\om_{(1)} (u_k (-1) \emu \ot E_{qk}(-1)) $$
$$ + \eru_{(n-3)} \sum_{q =1}^N m_q  
u_k (-1) \emu \ot \om_{(2)}  E_{qk}(-1))$$ 
$$ = \eru_{(n-1)} \sum_{q,k =1}^N m_q  D(u_k (-1) \emu \ot E_{qk}(-1))   
+ 2 \eru_{(n-2)} \sum_{q,k =1}^N m_q  u_k (-1) \emu \ot E_{qk}(-1)$$ 
$$- \eru_{(n-3)} D \emu.$$
The above yields zero when $n \geq 3$. Let us evaluate this expression
for $n = 0, 1, 2$. First we let $n = 2$:
$$\left( \om_{(-1)} \eru \right)_{(2)} 
\sum_{q,k =1}^N m_q  u_k (-1) \emu \ot E_{qk}(-1)
= - \eru_{(-1)} D \emu. $$
For $n = 1$ we get:
$$\left( \om_{(-1)} \eru \right)_{(1)} 
\sum_{q,k =1}^N m_q  u_k (-1) \emu \ot E_{qk}(-1)$$
$$=  2 \eru_{(-1)} \sum_{q,k =1}^N m_q  u_k (-1) \emu \ot E_{qk}(-1) 
- \eru_{(-2)} D \emu$$
$$=  2 \sum_{q,k =1}^N m_q  u_k (-1) \ermu \ot E_{qk}(-1) 
- (D\eru)_{(-1)} D \emu.$$
Finally, for $n = 0$ we obtain:
$$\left( \om_{(-1)} \eru \right)_{(0)} 
\sum_{q,k =1}^N m_q  u_k (-1) \emu \ot E_{qk}(-1)$$
$$ = \eru_{(-1)} D \left( \sum_{q,k =1}^N m_q  u_k (-1) \emu \ot E_{qk}(-1) \right)   
+ 2 (D\eru)_{(-1)} \sum_{q,k =1}^N m_q  u_k (-1) \emu \ot E_{qk}(-1) $$
$$- {1\over 2} (D^2 \eru)_{(-1)} D \emu.$$

{\it Step 3.}
The $n$-th products 
$$\left(
\sum_{p,s =1}^N r_p  u_s (-1) \eru \ot E_{ps}(-1) \right)_{(n)} 
\om_{(-1)} \emu $$
can be computed using the skew symmetry (\skeww) from the $n$-th
products computed in the previous step.
Indeed, letting 
$a = \sum\limits_{p,s =1}^N r_p  u_s (-1) \eru\ot E_{ps}(-1)$
and $b = \om_{(-1)} \emu$ and noticing that $b_{(j)} a = 0$ for
$j \geq 3$, we get from (\skeww):
$$ a_{(n)} b = 0 \quad \hbox{\rm for \ }n > 2,$$
$$ a_{(2)} b = - b_{(2)} a, $$ 
$$ a_{(1)} b = b_{(1)} a - D(b_{(2)} a), $$
$$ a_{(0)} b = - b_{(0)} a  + D(b_{(1)} a) 
- {1\over 2} D^2 (b_{(2)} a). $$
The expressions in the right hand sides are now available from Step 2
if we switch in those formulas $\r$ with $\m$. For $n = 2$ this will give us:
$$\left(
\sum_{p,s =1}^N r_p  u_s (-1) \eru \ot E_{ps}(-1) \right)_{(2)} 
\om_{(-1)} \emu  = \emu_{(-1)} D\eru .$$
For $n = 1$:
$$\left(
\sum_{p,s =1}^N r_p  u_s (-1) \eru \ot E_{ps}(-1) \right)_{(1)} 
\om_{(-1)} \emu  $$
$$= 2 \sum_{p,s =1}^N r_p  u_s (-1) \ermu \ot E_{ps}(-1) 
- (D\emu)_{(-1)} D \eru + D(\emu_{(-1)} D \eru) $$
$$= 2 \sum_{p,s =1}^N r_p  u_s (-1) \ermu \ot E_{ps}(-1) 
+ \emu_{(-1)} D^2 \eru .$$
And finally for $n = 0$:
$$\left(
\sum_{p,s =1}^N r_p  u_s (-1) \eru \ot E_{ps}(-1) \right)_{(0)} 
\om_{(-1)} \emu  $$
$$ = - \emu_{(-1)} D \left( 
\sum_{p,s =1}^N r_p  u_s (-1) \eru \ot E_{ps}(-1) \right)$$
$$ -2 (D\emu)_{(-1)} \sum_{p,s =1}^N r_p  u_s (-1) \eru \ot E_{ps}(-1)
+ {1 \over 2} (D^2 \emu)_{(-1)} D\eru$$
$$ + 2 D \left( \sum_{p,s =1}^N r_p  u_s (-1) \ermu \ot E_{ps}(-1) \right)
- D \left( (D \emu)_{(-1)} D\eru \right) + {1 \over 2} D^2 ( \emu_{(-1)} D\eru)$$ 
$$ = \emu_{(-1)} D\left( 
\sum_{p,s =1}^N r_p  u_s (-1) \eru \ot E_{ps}(-1) \right)
+ {1\over 2} \emu_{(-1)} D^3 \eru.$$

{\it Step 4.} The computation of the products
$$\left( \sum_{p,s =1}^N r_p 
 u_s (-1) \eru \ot E_{ps}(-1) \right)_{(n)} 
\sum_{q,k =1}^N m_q  u_k (-1) \emu \ot E_{qk}(-1) \eqno{(\BIV)} $$
will be based on Lemma \ltens (iii). Since we have chosen
the trivial cocycle on $\wgl$, $c_1 = c_I = 0$, then by Lemma \EE,
$E_{ps}(-1)_{(n)} E_{qk}(-1) = 0$ for $n \geq 1$. This implies that
(\BIV) is also zero for $n \geq 1$. The only case that we need to 
consider is $n=0$. Applying Lemma \ltens (iii) and Lemma \EE, we get
$$\left(
\sum_{p,s =1}^N r_p 
 u_s (-1) \eru \ot E_{ps}(-1) \right)_{(0)} 
\sum_{q,k =1}^N m_q u_k (-1) \emu \ot E_{qk}(-1) $$
$$ = \sum_{p,s,k =1}^N r_p m_s  u_s (-1) u_k (-1) \ermu \ot E_{pk}(-1)   
- \sum_{q,p,s =1}^N r_p m_q  u_s (-1) u_p (-1) \ermu \ot E_{qs}(-1)$$
$$ = (D\emu)_{(-1)} \sum_{p,k =1}^N r_p u_k(-1) \eru \ot E_{pk}(-1)
- (D\eru)_{(-1)} \sum_{q,s =1}^N m_q u_s(-1) \emu \ot E_{qs}(-1) .$$
This completes Step 4.

\

Now we can combine the results of the four steps and obtain the desired
$n$-th products (\dondo). For $n=3$ we get
$$\left( \om_{(-1)} \eru  + \sum_{p,s =1}^N r_p 
 u_s (-1) \eru \ot E_{ps}(-1) \right)_{(3)} 
\left( \om_{(-1)} \emu  + \sum_{q,k =1}^N m_q 
 u_k (-1) \emu \ot E_{qk}(-1) \right) $$
$$= 6 \ermu. $$
This corresponds via the commutator formula (\comm) to the term in 
(\ddo) with the third derivative of the delta function.

Collecting the terms with $n=2$ from all four steps we obtain
$$\left( \om_{(-1)} \eru  + \sum_{p,s =1}^N r_p 
 u_s (-1) \eru \ot E_{ps}(-1) \right)_{(2)} 
\left( \om_{(-1)} \emu  + \sum_{q,k =1}^N m_q 
 u_k (-1) \emu \ot E_{qk}(-1) \right) $$
$$ = 3 D\ermu - \eru_{(-1)} D\emu + (D\eru)_{(-1)} \emu =
4 (D\eru)_{(-1)} \emu + 2 \eru_{(-1)} D\emu .$$
These match the terms with the second derivative of the delta
function in (\ddo).

The case $n=1$ yields:
$$\left( \om_{(-1)} \eru  + \sum_{p,s =1}^N r_p 
 u_s (-1) \eru \ot E_{ps}(-1) \right)_{(1)} 
\left( \om_{(-1)} \emu  + \sum_{q,k =1}^N m_q 
 u_k (-1) \emu \ot E_{qk}(-1) \right) $$
$$ = 2 \om_{(-1)} \ermu + 2 (D\eru)_{(-1)} D\emu 
+ 2 \sum_{q,k =1}^N m_q  u_k (-1) \ermu \ot E_{qk}(-1) 
- (D\eru)_{(-1)} D \emu$$
$$+ 2 \sum_{p,s =1}^N r_p  u_s (-1) \ermu \ot E_{ps}(-1) 
+ \emu_{(-1)} D^2 \eru =$$
$$ = 2\left( \om_{(-1)} \ermu 
+ \sum_{p,s =1}^N (r_p + m_p)  u_s (-1) \ermu \ot E_{ps}(-1) 
\right) + D \left( (D\eru)_{(-1)}\emu \right) . $$
This is in agreement with the terms in (\ddo) with the first
derivative of the delta function.

Finally, for $n=0$:
$$\left( \om_{(-1)} \eru  + \sum_{p,s =1}^N r_p 
 u_s (-1) \eru \ot E_{ps}(-1) \right)_{(0)} 
\left( \om_{(-1)} \emu  + \sum_{q,k =1}^N m_q 
 u_k (-1) \emu \ot E_{qk}(-1) \right) $$
$$ = D(\om_{(-1)} \ermu) + {1\over 2} (D^2 \eru)_{(-1)} D\emu 
- {1\over 2} (D^3 \eru)_{(-1)} \emu $$
$$ + \eru_{(-1)} D \left( \sum_{q,k =1}^N m_q u_k (-1) \emu \ot E_{qk}(-1) \right)   
+ 2 (D\eru)_{(-1)} \sum_{q,k =1}^N m_q  u_k (-1) \emu \ot E_{qk}(-1) $$
$$- {1\over 2} (D^2 \eru)_{(-1)} D \emu 
 + \emu_{(-1)} D\left( 
\sum_{p,s =1}^N r_p  u_s (-1) \eru \ot E_{ps}(-1) \right)
+ {1\over 2} \emu_{(-1)} D^3 \eru $$
$$ + (D\emu)_{(-1)} \sum_{p,k =1}^N r_p u_k(-1) \eru \ot E_{pk}(-1)
- (D\eru)_{(-1)} \sum_{q,s =1}^N m_q u_s(-1) \emu \ot E_{qs}(-1) $$
$$ = D \left( \om_{(-1)} \ermu 
+ \sum_{p,s =1}^N (r_p + m_p)  u_s (-1) \ermu \ot E_{ps}(-1) 
\right) . $$
This matches the term in (\ddo) that comes with the delta function.

The proofs of both Lemma \Yddo \ and Theorem \main \ are now complete. 
 
\

{\bf 5. Irreducible representations of the full toroidal Lie algebra
and of its subalgebra with the divergence free vector fields.}

\

{\bf 5.1. Irreducible modules for the full toroidal Lie algebra.}

We are now going to apply the results of the previous section for the construction
of irreducible representations for the full toroidal Lie algebra $\g$.

\

{\bf Theorem \irre.} 
{\it
Let the constants $c, c_L, c_{LI}$ satisfy the assumptions
of Theorem \main. Let $L_\wdg(\lambda, c)$ be an irreducible highest weight
module for $\wdg$ with the highest weight $(\lambda, c)$, $\lambda \in \dh^*$. 
Let for $\alpha\in\C^N, \beta\in\Z^N$,
$M_\hyp^+ (\alpha,\beta)$ be the irreducible VOA module for $V_\hyp^+$, defined in (\MHyp).
Let $L_\wsl(\lambda_1, 0)$ be the irreducible highest weight
$\wsl$ module of level 0,
where $\lambda_1$ is a linear functional on the Cartan subalgebra of $\slN$.
 Let $L_\HVir(h,h_I,c_L,c_{LI},0)$ be the
irreducible highest weight module for the twisted Heisenberg-Virasoro
algebra. Then
$$L_{\tor} = L_\wdg(\lambda, c) \ot M_\hyp^+ (\alpha,\beta) \ot
L_\wsl(\lambda_1, 0) \ot L_\HVir(h,h_I,c_L,c_{LI},0) $$
has a structure of an irreducible module for the toroidal Lie algebra
$\g({1\over c},0)$.  
}

{\it Proof.} First of all we note that $L_{\tor}$ is an irreducible VOA module for
$$V_{\tor} = V_\wdg(c) \ot V_\hyp^+ \ot V_\wsl(0) \ot V_\HVir(c_L,c_{LI},0),$$
and thus by Corollary \tormod, $L_{\tor}$ is a module for the toroidal Lie
algebra $\g({1\over c},0)$. We are going to show that the fields (\ko)-(\tdd),
corresponding to $\g$, generate the VOA $V_{\tor}$. After this is done, we
use the Borcherds' identity (\qas) with $m=0$ to see that each moment of
every vertex operator from $V_{\tor}$ is a linear combination of associative
products of the operators corresponding to $\g$. Thus every $\g$ submodule is
also a submodule for $V_{\tor}$ and hence the simplicity of $L_{\tor}$
as a VOA module implies its irreducibility as a module for $\g$.

So we need to show that the set
$$S = \left\{ \eru, g(-1)\ot \eru, u_p(-1)\eru,
v_a(-1)\eru + \sum_{p=1}^N r_p \eru \ot E_{pa}(-1), \right.$$
$$\left. \om_{(-1)} \eru + \sum_{p,s=1}^N r_p u_s (-1)\eru \ot E_{ps}(-1) 
\right\}_{\r\in\Z^N} $$
generates $V_{\tor}$. Denote by $\S$  
the vertex subalgebra in $V_{\tor}$
generated by $S$. Since $V_{\tor}$ is a tensor product of four VOAs, it is
sufficient to show that each of the tensor factors, 
$V_\wdg, V_\hyp^+ , V_\wsl$ and $V_\HVir$, are in $\S$.

By taking $\r = 0$, we see that $g(-1), u_p(-1), v_p(-1)$ and $\om$ are in $S$.
Since $V_\wdg$ is spanned by $g_1(-n_1) \ldots g_k (-n_k) \o = 
g_1(-1)_{(-n_1)} \ldots g_k(-1)_{(-n_k)} \o$, we conclude that  
$V_\wdg \subset \S$. Since $V_\hyp^+$ is generated by $u_p(-1), v_p(-1), \eru$,
we get that $V_\hyp^+ \subset \S$. 
 To show that $V_\wsl \subset \S$, we choose $\r$ with $r_j = \delta_{pj}$, 
$j = 1,\ldots,N$. Then we get that $v_a(-1) e^{u_p} + e^{u_p}\ot E_{pa}(-1) \in S$.
Since $v_a(-1) e^{u_p} \in \S$, we obtain that $e^{u_p}\ot E_{pa}(-1) \in \S$
and $e^{-u_p}_{(-1)} e^{u_p}\ot E_{pa}(-1) = E_{pa}(-1) \in \S$. However
$V_\wsl$ is generated by $\{ E_{pa}(-1) | p\neq a \}$ and thus $V_\wsl \subset \S$.

For the last tensor factor we note that 
$\om = \om_{\wdg} + \om_\hyp + \om_\wsl + \om_\HVir \in S$ and also 
$\om_{\wdg}, \om_\hyp, \om_\wsl \in \S$, hence $\om_\HVir \in \S$.
Finally, we have $E_{aa} (-1) = \psi_1 (E_{aa})(-1) + \psi_2 (E_{aa})(-1) \in \S$,
where $\psi_1 (E_{aa})(-1) \in V_\wsl \subset \S$ and 
$\psi_2 (E_{aa})(-1) = {1 \over N} I(-1)$. But $I(-1)$ and $\om_\HVir$ generate
$V_\HVir$, thus $V_\HVir \subset \S$.

Since all four tensor factors are in $\S$, we conclude that $V_{\tor}$ is generated
by the set $S$. The theorem is now proved.

\

{\bf 5.2. Vertex operator algebra and irreducible representations for $\gdiv$.}

Now we will study the restriction of the modules for the toroidal Lie algebra to
the subalgebra $\gdiv = (\dg\otimes R) \oplus \K \oplus \Ddiv$. This subalgebra is
spanned by the elements $t_0^j \t^\r k_0$, $t_0^j \t^\r k_p$, $t_0^j \t^\r g, d_0$,
$t_0^j \t^\r d_a$ with $r_a = 0$ and
$$ t_0^j \t^\r \wda = - r_a t_0^j \t^\r d_0 + {1\over c} r_a (j + 1) t_0^j \t^\r k_0 
+ j t_0^j \t^\r d_a. \eqno{(\wdoa)}$$
The elements $t_0^j \t^\r k_0, t_0^j \t^\r k_p, t_0^j \t^\r g$ and 
$t_0^j \t^\r d_a$ with $r_a = 0$ correspond to the fields (\ko), (\ka),
(\ggg) and (\da) with $r_a = 0$. Just as in the previous theorem we get that the 
moments of these fields generate $V_\wdg, V_\hyp^+$ and $V_\wsl$. 
Collect the elements of the form (\wdoa) into the fields:
$$\wda (\r, z) = \sum_{j\in\Z} t_0^j \t^\r \wda z^{-j-2}.$$
Using (\tdd) and (\da) we obtain that $\wda (\r, z)$ is represented in the following way:
$$ \wda (\r, z) \mapsto r_a Y \left( \om_{(-1)} \eru + 
\sum_{p,s=1}^N r_p u_s(-1) \eru \ot E_{ps} (-1) , z \right)$$
$$- \left( z^{-1} + {\d \over \d z} \right) Y \left( v_a(-1) \eru + 
\sum_{p=1}^N r_p \eru \ot E_{pa} (-1) , z \right) 
+ r_a {z^{-2} \over 2} Y \left( \eru, z \right). \eqno{(\dwd)}$$

Since $V_\wdg \ot V_\hyp^+ \ot V_\wsl$ is generated by $\gdiv$, we shall consider
in (\dwd) only the terms that involve $V_\HVir$:
$$r_a Y \left( \eru \ot \om_\HVir, z \right) + {r_a \over N}
Y \left( \sum_{p=1}^N r_p u_p(-1) \eru \ot I(-1), z \right)$$
$$- {r_a \over N} Y \left( D(\eru \ot I(-1)), z \right)
- z^{-1} {r_a \over N} Y \left( \eru \ot I(-1), z \right)
+ r_a  {z^{-2} \over 2} Y \left( \eru, z \right)$$
$$ = r_a \left( Y \left( \eru \ot \om_\HVir, z \right)
- {1 \over N} Y \left( \eru \ot D(I(-1)), z \right) \right.$$
$$\left.
- z^{-1} {1 \over N} Y \left( \eru \ot I(-1), z \right)
+ {z^{-2} \over 2} Y \left( \eru, z \right) \right)$$
$$ = r_a Y \left(\eru, z \right) \ot
\left( Y \left(\om_\HVir, z \right)
- {1 \over N} Y \left( D(I(-1)), z \right)
- z^{-1} {1 \over N} Y \left( I(-1), z \right)
+ {z^{-2} \over 2} \Id \right) . \eqno{(\prj)}$$
To understand the structure of the expression 
$$Y \left(\om_\HVir, z \right)
- {1 \over N} \left(z^{-1} + {\d \over \d z} \right) 
Y \left( I(-1), z \right)
+ {z^{-2} \over 2} \Id, \eqno{(\expr)}$$
we consider the following

{\bf Lemma \emb.} 
{\it
Let $\wVir$ be the Virasoro algebra with the basis
$\left\{ \wL (n), \wCL \right\}$. For any $\gamma \in\C$ the map
$$\rho_\gamma : \quad \wVir \rightarrow \HVir / <C_{I}> ,$$
given by 
$$\rho_\gamma (\wL (n)) = L(n) + \gamma n I(n) + \delta_{n,0} \gamma C_{LI},$$
$$\rho_\gamma (\wCL) = C_L + 24 \gamma C_{LI} ,$$
is an embedding of Lie algebras.
}

\

{\bf Corollary \Memb.}
{\it
(a) Let $h, h_I, c_L, c_{LI} \in \C$. The homomorphism
$\rho_\gamma$ extends to the embedding of the Verma module
$M_\wVir (h + \gamma c_{LI}, c_L + 24 \gamma c_{LI})$ for the Virasoro algebra $\wVir$, 
into the Verma module 
$M_\HVir (h, h_I, c_L, c_{LI}, 0)$ for the twisted Heisenberg-Virasoro algebra $\HVir$.

\noindent
(b) Let $c_{LI} = {N \over 2}, \gamma = {1 \over N}$. Under the map $\rho_\gamma$
we have the correspondence of the fields 
$$\sum_{n\in\Z} \wL(n) z^{-n-2} \mapsto \sum_{n\in\Z} L(n) z^{-n-2}
- {1 \over N} \left( z^{-1} + {\d \over \d z} \right) \sum_{n\in\Z} I(n) z^{-n-1}
+ {z^{-2} \over 2} \Id . \eqno{(\embd)}$$
}

The proof of the Lemma is a straightforward computation, and the Corollary is an
immediate consequence. Note that the field in the right hand side in part (b) of
the Corollary coincides precisely with (\expr). Thus when $c_{LI} = {N \over 2}$
(as we have in Theorem \main), the components of this field satisfy the relations
of the Virasoro algebra with the value of the central charge $\wcL = c_L + 12$.

The field (\expr) involves vertex operators that are shifted by powers of $z$.
To deal with such expressions we need the following generalization of the
``preservation of identities'' principle (cf. [Li], Lemma 2.3.5):

{\bf Lemma \presr.}
{\it
Let $V$ be a VOA and let $M$ be a VOA module for $V$. Let $a^s, b^k, c^{nj} \in V$,
where $s,k,n,j$ run over finite subsets of $\Z$ and $n\geq 0$.

(i) If $M$ is a faithful VOA module and $\sum\limits_s z^{-s} Y_M (a^s, z) = 0$
(finite sum), then $a^s = 0$ for all $s$.

(ii) If 
$$\left[ \sum_s z^{-s} Y_V (a^s, z), \sum_k w^{-k} Y_V (b^k,w) \right]
= \sum_{n\geq 0} \sum_j w^{-j} Y_V (c^{nj},w) 
\left[ z^{-1} \left( {\d \over \d w} \right)^n \delta \left( {w\over z} \right) \right]$$
$$ \hbox{\hskip 4cm} \hbox{\it (all sums finite)} \eqno{(\relV)}$$
then
$$\left[ \sum_s z^{-s} Y_M (a^s, z), \sum_k w^{-k} Y_M (b^k,w) \right]
= \sum_{n\geq 0} \sum_j w^{-j} Y_M (c^{nj},w) 
\left[ z^{-1} \left( {\d \over \d w} \right)^n \delta \left( {w\over z} \right) \right].
\eqno{(\relM)}$$

(iii) If $M$ is a faithful VOA module then (\relM) implies (\relV).
}

{\it Proof.} Let us prove (i). Let $D_M$ be the infinitesimal translation operator on $M$.
If $\sum\limits_s Y_M (a^s, z) z^{-s} = 0$ then 
$$ 0 = z \sum\limits_s  \left[ D_M, Y_M (a^s, z)  \right] z^{-s} 
= z \sum_s \left( {\d \over \d z} Y_M (a^s, z) \right) z^{-s} $$
$$ = z \sum_s \left( {\d \over \d z} Y_M (a^s, z) \right) z^{-s}
- z {\d \over \d z} \left( \sum_s  Y_M (a^s, z) z^{-s} \right) $$
$$ = - \sum\limits_s Y_M (a^s, z) z {\d \over \d z} (z^{-s}) 
= \sum\limits_s s Y_M (a^s, z) z^{-s} .$$
Repeating this argument, we get that for any $m = 0,1,2, \ldots $
$$ \sum\limits_s s^m Y_M (a^s, z) z^{-s}  = 0 .$$
Since the sum in $s$ is finite, we can apply the Vandermonde determinant argument and
derive that $Y_M (a^s, z) = 0$ for all $s$. By the definition of the faithful module,
this implies that all $a^s = 0$.   

 To prove (ii), we use the commutator formula (\comm) and the basic properties of
the delta-function:
$$\left[ \sum_s z^{-s} Y_V (a^s, z), \sum_k w^{-k} Y_V (b^k,w) \right]$$
$$= \sum_{n,i\geq 0} \sum_{s,k} {1 \over n!} \pmatrix{ -s \cr i \cr}
 w^{-k-s-i} Y_V (a^s_{(n+i)} b^k, w) 
\left[ z^{-1} \left( {\d \over \d w} \right)^n \delta \left( {w\over z} \right) \right]
\quad \hbox{\rm (all sums finite)}. \eqno{(\relVY)}$$
By Corollary 2.2 from [K2], we obtain that for all $n \geq 0$,
$$ \sum_j w^{-j} Y_V (c^{nj},w) = 
\sum_{i\geq 0} \sum_{s,k} {1 \over n!} \pmatrix{ -s \cr i \cr}
 w^{-k-s-i} Y_V (a^s_{(n+i)} b^k, w).$$   
Since $V$ is a faithful VOA module over itself, we get using part (i) of the Lemma
that
$$ c^{nj} = \sum_{s,k} \sum_{{i\geq 0}\atop {s+k+i = j}} 
{1 \over n!} \pmatrix{ -s \cr i \cr} a^s_{(n+i)} b^k . \eqno{(\cab)}$$
However the relation (\relVY) holds in every VOA module M. Taking (\cab) into account,
we see that (\relM) also holds.

The proof for part (iii) is similar. We first see that
$$\sum_{n\geq 0} \sum_j w^{-j} Y_M (c^{nj},w) 
\left[ z^{-1} \left( {\d \over \d w} \right)^n \delta \left( {w\over z} \right) \right] =$$
$$\sum_{n,i\geq 0} \sum_{s,k} {1 \over n!} \pmatrix{ -s \cr i \cr}
 w^{-k-s-i} Y_M (a^s_{(n+i)} b^k, w) 
\left[ z^{-1} \left( {\d \over \d w} \right)^n \delta \left( {w\over z} \right) \right] .$$
Again using Corollary 2.2 from [K2] and part (i) of the Lemma, we obtain that the 
relation
(\cab) holds in $V$. Thus (\relV) also holds. This completes the proof of the Lemma.

\

Now we have done all the preparatory work and now ready to describe the representations 
for $\gdiv$.
  
{\bf Theorem \maind.}
{\it
Let $c, \wcL$ be complex numbers such that
$c\neq 0, c \neq -h^\vee$ and
$$ {c \dim \dg \over c + h^\vee} + 2N + \wcL = 24. $$
Let $M_\wdg$ be a VOA module for affine VOA $V_\wdg(c)$, 
$M_\hyp^+$ be a VOA module for the sub-VOA $V_\hyp^+$
of the hyperbolic lattice VOA, 
$M_\wsl$ be a module for affine $\wsl$ VOA $V_\wsl(0)$,
and $M_{\wVir}$ be a VOA module for the Virasoro VOA
$V_{\wVir}(\wcL)$. Then 
$$M_\gdiv = M_\wdg \ot M_\hyp^+ \ot M_{\wsl} \ot M_{\wVir}$$
is a module for the divergence free subalgebra $\gdiv({1\over c})$ 
of the toroidal Lie algebra.
The action of the fields
$k_0 (\r,z), k_a(\r,z), g(\r,z)$ is given by
the formulas (\ko),(\ka),(\ggg).
The action of $d_a(\r,z)$ with $r_a = 0$ is given by (\da).
The action of $d_0$ is given by 
$$ d_0 \mapsto \Id - \om_{(1)}, \eqno{(\doo)}$$
where $\om$ is the Virasoro element of the tensor product VOA
$$V_\gdiv = V_\wdg(c) \ot V_\hyp^+ \ot V_\wsl(0) \ot V_{\wVir}(\wcL).$$
Finally, the field $\wda(\r,z)$ is represented by 
$$ r_a Y_M \left(\om_{(-1)} \eru + \sum_{p,s =1}^N r_p u_s (-1) \eru \ot 
\psi_1 (E_{ps}) (-1), z \right) $$
$$ - \left( z^{-1} + {\d \over \d z} \right)
Y_M \left( v_a(-1) \eru + \sum_{p=1}^N r_p \eru \ot 
\psi_1 (E_{pa}) (-1), z \right).  \eqno{(\wdaa)}$$
}
  
{\it Proof.} The Lie bracket in $\gdiv$ may be encoded in the commutator
relations between the fields $k_0 (\r,z), k_a(\r,z), g(\r,z)$, $\wda(\r,z)$
$d_a(\r,z)$ with $r_a = 0$, and the element $d_0$, analogous to (\dko)-(\ddo).
We need to show that the same commutator relations hold for their images
(\ko), (\ka), (\ggg), (\da) with $r_a = 0$, (\doo) and (\wdaa). 
It is easy to see that the relations involving $d_0$ in the left hand sides, 
hold due to $(\omdeg)$.
Also, $d_0$ does not belong to the commutant of $\gdiv$ and will not appear
in the right hand sides of the commutator relations.
The commutator relations that should be verified for the remaining fields
(\ko), (\ka), (\ggg), (\da) with $r_a = 0$, and (\wdaa)
are of the form (\relM). 
Our strategy is to embed one of the modules for $V_\gdiv$ into a module for
the full toroidal Lie algebra $\g$. This embedding will have the property
that the restriction of the action of $\g$ to subalgebra $\gdiv$ will coincide
with (\ko), (\ka), (\ggg), (\da) with $r_a = 0$, (\doo) and (\wdaa).
This will imply that the necessary commutator relations hold in the chosen
module for $V_\gdiv$. Since the module that we will consider will be faithful,
then by preservation of identities (Lemma \presr), the same required relations
will hold in all VOA modules for $V_\gdiv$.

Let us carry out this plan. 
Consider the embedding  given by Corollary \Memb \ with
$h = -{1 \over 2}$, $\gamma = {1 \over N}, c_{LI} = {N \over 2}$, $h_I = c_I =0$,
of the Verma module 
$M_\wVir (0,\wcL)$ for the Lie algebra $\wVir$, 
into the Verma module $M_\HVir (-{1\over 2}, 0, \wcL - 12, {N \over 2}, 0)$
for the twisted Heisenberg-Virasoro algebra. Under this homomorphism
$$ \wL (z) = \sum_{n\in\Z} \wL (n) z^{-n-2} \mapsto
L(z) - {1 \over N} \left( z^{-1} + {\d \over \d z} \right) I(z) + 
{z^{-2} \over 2} \Id.$$
This map extends to the embedding 
$$ V_\wdg(c) \ot V_\hyp^+ \ot V_\wsl(0) \ot M_{\wVir}(0,\wcL) \subset
V_\wdg(c) \ot V_\hyp^+ \ot V_\wsl(0) \ot 
M_\HVir (-{1\over 2}, 0, \wcL - 12, {N \over 2}, 0).$$
By Corollary \tormod, the latter is a module for the full toroidal Lie algebra 
$\g({1\over c}, 0)$. We consider the restriction of this representation
to the subalgebra $\gdiv({1\over c})$ and claim that 
$ V_\wdg(c) \ot V_\hyp^+ \ot V_\wsl(0) \ot M_{\wVir}(0,\wcL)$ is invariant
under the action of $\gdiv({1\over c})$. The action of 
$k_0 (\r,z)$, $k_a(\r,z)$, $g(\r,z)$ and $d_a(\r,z)$ with $r_a = 0$ is given by
(\ko),(\ka),(\ggg) and (\da). 
Let us show that the action of $\wda(\r,z)$ on 
$V_\wdg(c) \ot V_\hyp^+ \ot V_\wsl(0) \ot 
M_\HVir (-{1\over 2}, 0, \wcL - 12, {N \over 2}, 0)$
coincides with (\wdaa).
Indeed, following the computations (\dwd)--(\embd), we get:
$$\wda(\r,z) \mapsto r_a Y \left(
( \om_\wdg + \om_\Hyp + \om_\wsl + \om_{\HVir} )_{(-1)} \eru
+ \sum_{p,s = 1}^N r_p u_s(-1) \eru \ot E_{ps} (-1), z \right) $$
$$ - \left( z^{-1} + {\d \over \d z} \right) 
Y \left( v_a (-1) \eru + \sum_{p=1}^N r_p \eru \ot E_{pa} (-1), z \right)
+ r_a {z^{-2} \over 2} Y(\eru, z)$$
$$ = r_a Y \left(
( \om_\wdg + \om_\Hyp + \om_\wsl)_{(-1)} \eru
+ \sum_{p,s = 1}^N r_p u_s(-1) \eru \ot \psi_1(E_{ps}) (-1), z \right) $$
$$ - \left( z^{-1} + {\d \over \d z} \right) 
Y \left( v_a (-1) \eru + \sum_{p=1}^N r_p \eru \ot \psi_1 (E_{pa}) (-1), z \right)$$
$$ + r_a Y(\eru \ot \om_\HVir , z) + r_a {1 \over N} Y( (D\eru) \ot I(-1), z)$$
$$ - r_a \left( z^{-1} + {\d \over \d z} \right) {1 \over N}
Y (\eru \ot I(-1) , z) + r_a {z^{-2} \over 2} Y(\eru, z)$$
$$ = r_a Y \left(
( \om_\wdg + \om_\Hyp + \om_\wsl + \om_{\wVir} )_{(-1)} \eru
+ \sum_{p,s = 1}^N r_p u_s(-1) \eru \ot \psi_1 (E_{ps}) (-1), z \right) $$
$$ - \left( z^{-1} + {\d \over \d z} \right) 
Y \left( v_a (-1) \eru + \sum_{p=1}^N r_p \eru \ot \psi_1 (E_{pa})(-1), z \right),$$
which is the same as (\wdaa).
Thus the specified action defines a 
representation of $\gdiv$ on $V_\wdg(c) \ot V_\hyp^+ \ot V_\wsl(0) \ot M_{\wVir}(0,\wcL)$,
and the fields (\ko), (\ka), (\ggg), (\da) with $r_a = 0$, and (\wdaa) satisfy the
relations that reflect the Lie bracket in $\gdiv$.
This module is a faithful VOA module for $V_\gdiv$, since $V_\gdiv$ itself is its
factor module. Thus by the preservation of identities, Lemma \presr,
the required commutator
relations hold in $V_\gdiv$ and in all VOA modules for $V_\gdiv$. This establishes
the claim of the theorem.

\

In the next theorem we give the description of the irreducible modules for $\gdiv$.

\

{\bf Theorem \irrd.} 
{\it
Let constants $c, \wcL$ satisfy the assumptions of Theorem
\maind. 
Let $L_\wdg(\lambda, c)$ be an irreducible highest weight
module for $\wdg$. Let for $\alpha\in\C^N, \beta\in\Z^N$,
$M_\hyp^+ (\alpha,\beta)$ be the irreducible VOA module for $V_\hyp^+$, defined in (\MHyp).
Let $L_\wsl(\lambda_1, 0)$ be the irreducible highest weight
$\wsl$ module of level 0,
where $\lambda_1$ is a linear functional on the Cartan subalgebra of $\slN$.
 Let $L_\wVir(h,\wcL)$ be the
irreducible highest weight module for the Virasoro Lie algebra, $h\in\C$. 
Then
$$L_{\gdiv} = L_\wdg(\lambda, c) \ot M_\hyp^+ (\alpha,\beta) \ot
L_\wsl(\lambda_1, 0) \ot L_\wVir(h,\wcL) $$
has a structure of an irreducible module for the Lie algebra
$\gdiv({1\over c})$.  
}

The proof of this theorem is completely analogous to the proof of 
Theorem \irre \ and will be omitted.

\

We conclude the paper with two observations. The relation for the central charges
for the modules for the $N+1$-toroidal Lie algebra with the divergence-free
vector fields $\gdiv$ constructed in Theorem \maind, may be rewritten as
$$ {c \dim \dg \over c + h^\vee} + 2(N+1) + \wcL = 26. \eqno{(\strng)}$$
This has a striking resemblance to the ``no ghost'' theorem in string theory.
If we choose the modules for the affine algebras $\wdg$, $\wsl$ and for the Virasoro
algebra to be unitary, this would require $c > 0$ and $\wcL \geq 0$ and thus
we get that such modules exist only when
$$ N < 12.$$
For example, if we choose the basic module for $\wdg$ and trivial modules for
$\wsl$ and for the Virasoro algebra then the condition (\strng) will become
$$ \rank (\dg) + 2(N+1) = 26 .$$
Unfortunately, even when the affine and the Virasoro parts are unitary, the module
for $\gdiv$ is not unitarizable because the hyperbolic lattice is not positive-definite
and the lattice VOA used for the construction of such modules does not possess a unitary
structure.

 The second curious fact is that at the critical value $N = 12$ we get an exceptional
module for the Lie algebra
$$ \Ddiv \oplus \K.$$
The Lie bracket in this algebra is given by (\Ldk) and (\Ldd) with cocycle $\tau_1$.
Note that this Lie algebra has a non-degenerate symmetric invariant bilinear form
given by (\pairDK).
If we take trivial modules for the Virasoro and for the affine algebras $\wdg$, $\wsl$,
then only when $N=12$ we arrive at the representation of $ \Ddiv \oplus \K$ just
on the lattice part $V^+_\hyp$.
We get the following remarkable result:

{\bf Theorem \DK .} 
{\it
Let $N = 12.$ Then $V^+_\hyp$ has a structure of a module for
$\Ddiv \oplus \K$ with cocycle $\tau_1$. 
The action of the Lie algebra is given by
$$ k_0 (\r, z) \mapsto Y(\eru, z), \quad
k_p (\r, z) \mapsto Y( u_p(-1) \eru, z), $$
$$ d_0 \mapsto \Id - \om_{(1)} , \quad d_a \mapsto v_a (0), $$
$$\wda (\r, z) \mapsto r_a Y(\om_{(-1)} \eru, z) - 
\left( z^{-1} + {\d \over \d z} \right) Y (v_a(-1) \eru, z).$$
}

The character of this module with respect to the diagonalizable operators 
$d_0$, $d_1$, $\ldots$, $d_N$ has nice modular properties -- 
it is a product of $12$ delta-functions with the 
$-24$-th power of the Dedekind $\eta$-function:
$$ \char V^+_\hyp = q_0 \prod_{k=1}^\infty \left( 1 - q_0^{-k} \right)^{-24} 
\times \prod_{p=1}^{12} \sum_{j \in\Z} q_p^j .$$ 
 
\

\

{\bf References:}

\item{[AABGP]}  Allison, B.N., Azam, S., Berman, S., Gao, Y., Pianzola, A.: 
{\it Extended affine Lie algebras and their root systems.} 
Mem.Amer.Math.Soc. {\bf 126}, no. 603, 1997.

\item{[ACKP]} Arbarello, E., De Concini, C., Kac, V.G., Procesi, C.:
{Moduli spaces of curves and representation theory.}
Commun.Math.Phys. {\bf 117}, 1-36 (1988).

\item{[BB]} Berman, S., Billig, Y.:
{Irreducible representations for toroidal Lie algebras.} 
J.Algebra {\bf 221}, 188-231 (1999).

\item{[BBS]} Berman, S., Billig, Y., Szmigielski, J.:
{Vertex operator algebras and the representation theory of toroidal algebras.}
to appear in Contemp.Math.

\item{[BC]} Berman, S., Cox, B.: 
{ Enveloping algebras and representations of toroidal Lie algebras.} 
Pacific J.Math. {\bf 165}, 239-267 (1994).

\item {[BGK]} Berman, S., Gao, Y., Krylyuk, Y.: 
{Quantum tori and the structure of elliptic quasi-simple Lie algebras.} 
J.Funct.Analysis {\bf 135}, 339-389 (1996).

\item{[B1]} Billig, Y.: 
{Principal vertex operator representations for toroidal Lie algebras.}  
J. Math. Phys. {\bf 39}, 3844-3864 (1998).

\item{[B2]} Billig, Y.: 
{An extension of the KdV hierarchy arising from a representation 
of a toroidal Lie algebra.} 
J.Algebra {\bf 217}, 40-64 (1999).

\item{[B3]} Billig, Y.:
{Representations of the twisted Heisenberg-Virasoro algebra at level 
zero.}
preprint.

\item{[DLM]} Dong, C., Li, H., Mason, G.:
{Vertex Lie algebras, vertex Poisson algebras and vertex algebras.}
preprint, QA-0102127.

\item{[EM]} Eswara Rao, S., Moody, R.V.: 
{Vertex representations for $n$-toroidal
Lie algebras and a generalization of the Virasoro algebra.}
Commun.Math.Phys. {\bf 159}, 239-264 (1994).

\item{[FKRW]} Frenkel, E., Kac, V., Radul, A., Wang, W.:
{$W_{1+\infty}$ and $W(gl_\infty)$ with central charge $N$.}
Commun.Math.Phys. {\bf 170}, 337-357 (1995).

\item{[FLM]} Frenkel, I.B., Lepowsky, J., Meurman, A.:
{\it Vertex operator algebras and the Monster.} 
New York: Academic Press, 1988.

\item{[IT]} Ikeda, T., Takasaki, K.:
{Toroidal Lie algebras and Bogoyavlensky's 2+1-dimensional equation.}
Internat.Math.Res.Notices {\bf 7}, 329-369 (2001).

\item{[IKU]} Inami, T., Kanno, H., Ueno, T.:
{Higher-dimensional WZW model on K\"ahler manifold and toroidal Lie algebra.}
Mod.Phys.Lett. A {\bf 12}, 2757-2764 (1997).

\item{[IKUX]} Inami, T., Kanno, H., Ueno, T., Xiong, C.-S.:
{Two-toroidal Lie algebra as current algebra of four-dimensional K\"ahler WZW model.}
Phys.Lett. B {\bf 399}, 97-104 (1997).

\item{[ISW]} Iohara, K., Saito, Y., Wakimoto, M.: 
{Hirota bilinear forms with 2-toroidal symmetry.} 
Phys.Lett. A {\bf 254}, 37-46 (1999).

\item{[K1]} Kac, V.: 
{\it Infinite dimensional Lie algebras.} 
Cambridge: Cambridge University Press, 3rd edition, 1990.

\item{[K2]}  Kac, V.: 
{\it Vertex algebras for beginners.} 
Second Edition, University Lecture Series, {\bf 10}, A.M.S., 1998.

\item{[Kas]} Kassel, C.: 
{K\"ahler differentials and coverings of complex simple 
Lie algebras extended over a commutative ring.}
J.Pure Applied Algebra {\bf 34}, 265-275 (1984).

\item{[L]} Larsson, T.A.:
{Lowest-energy representations of non-centrally extended diffeomorphism 
algebras.}
Commun.Math.Phys. {\bf 201}, 461-470 (1999).

\item{[Li]} Li, H.:
{Local systems of vertex operators, vertex superalgebras and modules.}
{J.Pure Appl.Algebra} {\bf 109}, 143-195 (1996).

\item{[MRY]} Moody, R.V., Rao, S.E., Yokonuma, T.:  
{Toroidal Lie algebras and vertex representations.} 
Geom.Ded. {\bf 35}, 283-307 (1990).

\item{[P]} Primc, M.:
{Vertex algebras generated by Lie algebras.}
J.Pure Appl.Algebra {\bf 135}, 253-293 (1999).

\item{[R]} Roitman, R.:
{On free conformal and vertex algebras.}
J.Algebra {\bf 217}, 496-527 (1999).

\end